 \documentclass[12pt]{article}

\usepackage{graphicx}
\usepackage{latexsym,amsmath,amsfonts,amscd, amsthm}
\usepackage{changebar}
\usepackage{color}
\usepackage{bm}
\usepackage{tikz}
\usepackage{multirow}
\usepackage{epsfig,subfigure}
\usepackage{multirow}
\usetikzlibrary{arrows,backgrounds,snakes}

\newtheoremstyle{plainNoItalics}{}{}{\normalfont}{}{\bfseries}{.}{ }{}

\theoremstyle{plain}
\newtheorem{thm}{Theorem}[section]

\theoremstyle{plainNoItalics}

\newtheorem{defn}[thm]{Definition}

\newtheorem{prop}[thm]{Proposition}
\newtheorem{exa}[thm]{Example}

\newcommand{\beq}{\begin{equation}}
\newcommand{\eeq}{\end{equation}}
\newcommand{\beqa}{\begin{eqnarray}}
\newcommand{\eeqa}{\end{eqnarray}}
\newcommand{\bit}{\begin{itemize}}
\newcommand{\eit}{\end{itemize}}
\newcommand{\bedef}{\begin{defn}}
\newcommand{\edefn}{\end{defn}}
\newcommand{\bpro}{\begin{prop}}
\newcommand{\epro}{\end{prop}}

\begin{document}

\baselineskip=1.6pc

%\vspace*{.10in}

%=============  title  =========================
\begin{center}
{\bf
A Maximum-Principle-Satisfying High-order
Finite Volume Compact WENO Scheme for Scalar Conservation
Laws}

%with Applications in Incompressible
%Flow
\end{center}
\vspace{.2in}
\centerline{
Yan Guo \footnote{Department of Mathematics, China University of Mining and Technology, Xuzhou, Jiangsu 221116, P.R. China. Email: yanguo@cumt.edu.cn}
Tao Xiong \footnote{Department of
Mathematics, University of Houston, Houston, 77004, USA. E-mail:
txiong@math.uh.edu}
Yufeng Shi \footnote{School of Electric Power Engineering, China University of Mining and Technology, Xuzhou, Jiangsu 221116, P.R. China. Email: shiyufeng@cumt.edu.cn}
}

\bigskip
\centerline{\bf Abstract}
\vspace{.1in}
{In this paper, a maximum-principle-satisfying
finite volume compact scheme is proposed for solving scalar hyperbolic conservation laws. The scheme combines WENO schemes (Weighted Essentially Non-Oscillatory) with a class of compact schemes under a finite volume framework, in which the nonlinear WENO weights are coupled with lower order compact stencils. The maximum-principle-satisfying polynomial rescaling limiter in \cite{zhang2010maximum} is adopted to construct the present schemes at each stage of an explicit
Runge-Kutta method, without destroying high order accuracy and conservativity. Numerical examples
for one and two dimensional problems including incompressible flows are presented to assess the good
performance, maximum principle preserving, essentially non-oscillatory and highly accurate resolution
of the proposed method.
}

%It is well known that conservative compact finite
%volume schemes have high resolution properties and WENO (Weighted
%Essentially Non-Oscillatory) schemes are essentially oscillation
%free near the flow discontinuities. We extend the main idea of the
%WENO schemes to a class of classical compact finite volume schemes,
%where lower order compact stencils are combined with WENO nonlinear
%weights to get a higher order finite volume compact-WENO scheme. The
%maximum-principle-satisfying limiter is adopted to construct the
%present high order finite volume schemes in each stage for a
%Runge-Kutta method without destroying accuracy and conservativity. A
%number of test cases are presented to demonstrate the high-order
%accuracy, robustness, essentially non-oscillatory and
%high-resolution properties of the proposed scheme.}
\vfill

\newpage

\newpage

\section{Introduction}
\label{intro}
\setcounter{equation}{0}
\setcounter{figure}{0}
\setcounter{table}{0}

In this paper we consider the following scalar
hyperbolic conservation equation
\begin{equation}\label{adveqn}
\frac{\partial{u(x,t)}}{\partial t}+\nabla \cdot f(u(x,t))=0, \quad
x=(x_1,x_2,...,x_d) \in \mathbb{R}^d,
\end{equation}
with the initial condition $u(0,x)=u_0(x)$.
A main property of \eqref{adveqn} is that the solution $u(x,t)$ might
develop discontinuities in finite time even when the
initial data is smooth, due to its different propagation speeds. The weak solutions
of \eqref{adveqn} are not unique, hence a physically relevant
entropy solution should be considered. An important property of
the entropy solution for \eqref{adveqn} is that it satisfies a
strict maximum-principle \cite{dafermos2000conservation}, namely
\begin{equation}
u_m\leq u(x,t)\leq u_M, \quad \textrm{if $u_m\leq u(x,0)\leq u_M$},
\end{equation}
where $u_m=\min_x u(x,0)$ and $u_M=\max_x u(x,0)$. The total
variation diminishing (TVD) schemes \cite{harten1983high} satisfy
the strict maximum-principle, but they degenerate to first order
accuracy at smooth extrema \cite{osher1984high}. Recently Zhang et. al.
proposed uniformly high order accurate schemes satisfying a strict
maximum-principle based on the finite volume weighted essentially non-oscillatory (WENO) and finite element discontinuous Galerkin (DG) frameworks
for scalar conservation laws \cite{zhang2010maximum,zhang2011maximum}. These high
order schemes achieve the strict maximum principle by applying a polynomial rescaling
limiter at each stage of an explicit Runge-Kutta (RK) method or at each step of a multistep
method. The technique was later generalized to positivity
preserving high order DG and finite volume WENO schemes for compressible Euler equations \cite{zhang2010positivity}. Another class of high order parametrized maximum-principle-preserving (MPP)
flux limiters was proposed by Xu et. al. \cite{mpp_xu,mpp_xuMD} under a finite volume framework, which limits a high order numerical flux towards a first order monotone flux. Later in \cite{xiong2013parametrized} Xiong et. al.
generalized the parametrized high order MPP flux limiters for finite difference RK-WENO schemes with applications in incompressible flows. They only applied the limiters at the final stage of an explicit RK method which could save much computational cost. For these finite difference and finite volume WENO schemes solving
scalar hyperbolic conservation laws, they are high order
accurate in smooth regions and essentially non-oscillatory
for shock capturing. However, these methods based on non-compact WENO
schemes often suffer from excessive numerical dissipation, poor
spectral resolution and increasingly wide stencils with increasing order of accuracy.

A class of compact finite difference schemes was proposed
in \cite{lele1992compact}, which have significant higher spectral resolutions
with narrower stencils.
%and easier application of boundary conditions.
The compact scheme has been applied to incompressible flows
in \cite{liu1996essentially, wilson1998higher} and compressible flows in
\cite{lerat2001residual, ekaterinaris1999implicit}. However this classical linear compact finite difference
scheme often yields oscillatory solutions across discontinuities. To address this difficulty,
several hybrid schemes are proposed to couple the ENO and WENO schemes for simulating shock-turbulence interaction problems, e.g., a hybrid compact-ENO scheme by Adams et. al. \cite{adams1996high} and a hybrid compact WENO scheme by Pirozzoli \cite{pirozzoli2002conservative}.
Ren et. al. proposed a characteristic-wise hybrid compact WENO scheme \cite{ren2003characteristic} as a weighted average of the conservative compact scheme
\cite{pirozzoli2002conservative} and the WENO scheme \cite{jiang1996efficient}. These hybrid schemes used  smooth indicators to transit from the compact scheme to the WENO scheme around the discontinuities.

As an alternative to the hybrid schemes, a nonlinear compact scheme was proposed by Cockburn and Shu \cite{cockburn1994nonlinearly} based on TVD and TVB limiters. Deng et. al. then developed fourth order and fifth order weighted compact nonlinear schemes from cell-centered compact schemes and compact  interpolations of conservative variables at cell edges \cite{deng1997compact, deng2000developing}. The idea was generalized to weighted compact nonlinear schemes with increasingly high order accuracies in  \cite{zhang2008development} by directly interpolating the flux. However, these schemes are not truly compact schemes and the spectral resolution would be reduced. A new class of linear central compact scheme with spectral-like resolution was recently proposed in \cite{liu2013new} based on the cell-centered compact scheme of Lele \cite{lele1992compact}. Instead of interpolating the values on cell centers, they
directly evolve the values both on the grid nodes and the cell centers.

Another type of weighted compact nonlinear scheme was constructed in \cite{jiang2001weighted}. The scheme
was a weighted combination of compact substencils which were two biased third order
compact stencils and a central fourth order compact stencil. The scheme would result in a sixth order central
compact scheme if with optimal weights for smooth solutions. Ghosh and Baeder employed the idea and developed a new class of compact reconstruction finite difference WENO schemes \cite{ghosh2012compact}. In their approach,
lower order biased compact stencils were used and the higher order interpolation with optimal weights was
upwind. A tridiagonal system was solved at each time step. However, the scheme was high order accurate,
essentially non-oscillatory around discontinuities and superior spectral accurate.

Most of the above mentioned compact schemes are based on a finite difference framework.
Compact schemes based on a finite volume framework would be more nature especially
for unstructured meshes. In \cite{gaitonde1997optimized}, Gaitonde et. al. proposed
compact-difference-based finite-volume schemes for linear
wave phenomena. Kobayashi extended the work \cite{gaitonde1997optimized} to
a class of Pad\'{e} finite volume methods \cite{kobayashi1999class}.
A fourth-order finite volume compact scheme was provided for the incompressible
Navier-Stokes equations in \cite{pereira2001fourth} and applied to
incompressible Navier-Stokes equations on staggered grids \cite{hokpunna2010compact}
and compressible Navier-Stokes equations on nonuniform grids \cite{ghadimi2012fourth}.
Piller and Stalio developed compact finite volume schemes for one and two dimensional
transport and Navier-Stokes equations on stagger grids \cite{piller2004finite} and for
three dimensional scalar advection-diffusion equation on boundary-fitted grids \cite{piller2008compact},
among many others.

In this paper, we follow \cite{ghosh2012compact} to develop a finite volume compact WENO (FVCW) scheme,
where lower order compact stencils based on cell averages are used. The scheme combines the nonlinear
WENO weights to yield a fifth order upwind compact finite volume interpolation. As a new ingredient, we
incorporate the MPP polynomial rescaling limiter \cite{zhang2010maximum,zhang2011maximum},
which would be essentially important for some extreme problems with complex structures. A similar idea with positivity preserving (PP) limiter \cite{zhang2011maximum} for one dimensional compressible Euler
system has been explored in \cite{guo2014positivity}. For the FVCW scheme with the MPP limiter, numerical experiments will be presented for one and two dimensional scalar hyperbolic equations with application to incompressible Euler equations. The numerical results will show the good performance, maximum principle preserving and high resolution property of our new proposed approach.

%In this paper, a maximum-principle satisfying finite volume
%compact WENO (FVCW) scheme is constructed for two dimensional
%incompressible Flow. We employ the main idea that is described in
%\cite{guo2014positivity} for compressible Euler equations where
%lower order compact stencils are combined with WENO weights to yield
%a fifth-order upwind compact interpolation, and
%positivity-preserving limiter is used to preserve positive density
%and internal energy for compressible Euler equations. As an
%alternative to the finite difference compact schemes proposed in
%\cite{ghosh2012compact}, we developed a maximum-principle satisfying
%FVCW scheme where a linear scaling limiter \cite{zhang2011maximum}
%is added to the fifth-order finite volume compact WENO scheme to
%achieve strict maximum principle preserving. Comparing to classical
%small length scale fifth order finite volume compact schemes, the
%present schemes avoid spurious numerical oscillations in
%discontinuous region. The present small length scale FVCW scheme
%could also produce superior resolution and lower truncation errors
%compared to the classical maximum-principle satisfying finite volume
%WENO schemes.

The rest of the paper is organized as follows. In Section 2, the finite volume
compact WENO scheme is presented for the scalar conservation law in
one dimension and the maximum principle satisfying limiter is
introduced. The scheme for two dimensional conservation law on rectangular meshes
will be described in Section 3 and the application to incompressible flows
will be discussed in Section 4. In Section 5, we will show the numerical results
for scalar conservation problems and incompressible Euler
equations. Finally the conclusions are made in Section 6.

\section{Maximum-principle-satisfying finite volume compact WENO schemes}
\label{sec2}
\setcounter{equation}{0}
\setcounter{figure}{0}
\setcounter{table}{0}
\subsection{Finite volume compact WENO scheme}
We first briefly review the finite volume compact WENO scheme
\cite{guo2014positivity} for one dimensional hyperbolic conservation
equation
\begin{equation}\label{1Dadveqn}
\frac{\partial{u}}{\partial t}+\frac{\partial{f(u)}}{\partial x}=0,
\end{equation}
with the initial condition $u(0,x)=u_0(x)$. The computational domain
$[a,b]$ is divided into $N$ cells
\begin{equation*}
a=x_{\frac{1}{2}}<x_{\frac{3}{2}}<\cdots<x_{N+\frac{1}{2}}=b.
\end{equation*}
The cell is denoted by
$I_j=[x_{j-\frac{1}{2}},x_{j+\frac{1}{2}}]$. For simplicity, we consider an
uniform grid and the size of the cell is $\Delta x=\frac{b-a}{N}$.
The cell-averaged value of cell $j$, which denotes to be $\bar{u}_j$, can be defined as
\begin{equation}
\bar{u}(x_j,t)=\frac{1}{\Delta
{x}}\int^{x_{j+\frac{1}{2}}}_{x_{j-\frac{1}{2}}}u(x,t)dx.
\end{equation}
We approximate \eqref{1Dadveqn} by the following finite volume
conservative scheme
\begin{equation}\label{1DadvFV-b}
\frac{d\bar{u}_j(t)}{dt}=-\frac{1}{\Delta
{x}}(\hat{f}_{j+\frac{1}{2}}-\hat{f}_{j-\frac{1}{2}}),
\end{equation}
where the numerical flux $\hat{f}_{j+\frac{1}{2}}$ is defined by
\begin{equation}\label{flux}
\hat{f}_{j+\frac{1}{2}}=h(u^{-}_{j+\frac{1}{2}},u^{+}_{j+\frac{1}{2}}).
\end{equation}
Let $u_{j+\frac{1}{2}}^{-}$ denote a fifth order approximation of
the nodal value $u(x_{j+\frac{1}{2}},t^n)$ in cell $I_j$ and
$u_{j+\frac{1}{2}}^{+}$ denote a fifth order approximation of the
nodal value $u(x_{j+\frac{1}{2}},t^n)$ from cell $I_{j+1}$. In this
paper, $u^{-}_{j+\frac{1}{2}}$ and $u^{+}_{j+\frac{1}{2}}$ are
obtained from a high order compact WENO reconstruction, which will
be discussed in the following.

An optimal fifth-order compact upwind scheme \cite{pirozzoli2002conservative}
can be written as
\begin{equation}\label{5th-order1}
 \frac{3}{10}u_{j-\frac{1}{2}}^{-}+\frac{6}{10}u_{j+\frac{1}{2}}^{-}+\frac{1}{10}u_{j+\frac{3}{2}}^{-}
 =\frac{1}{30}\bar{u}_{j-1}+\frac{19}{30}\bar{u}_{j}+\frac{10}{30}\bar{u}_{j+1}.
\end{equation}
The classical small length scale finite volume linear compact scheme \eqref{5th-order1}
is very accurate and keep good resolutions in smooth regions. However, nonphysical
oscillations are generated when they are directly applied to problems with discontinuities
and the amplitude does not decrease as refining the grid. \eqref{5th-order1} is a
weighted combination of three third order compact substencils, they are
\begin{subequations}
\label{substencil}
\begin{align}
\frac{2}{3}{u}_{j-\frac{1}{2}}^-+\frac{1}{3}{u}_{j+\frac{1}{2}}^-&=\frac{1}{6}(\bar{u}_{j-1}+5\bar{u}_{j}), \\
\frac{1}{3}{u}_{j-\frac{1}{2}}^-+\frac{2}{3}{u}_{j+\frac{1}{2}}^-&=\frac{1}{6}(5\bar{u}_{j}+\bar{u}_{j+1}), \\
\frac{2}{3}{u}_{j+\frac{1}{2}}^-+\frac{1}{3}{u}_{j+\frac{3}{2}}^-&=\frac{1}{6}(\bar{u}_{j}+5\bar{u}_{j+1}).
\end{align}
\end{subequations}
and to get \eqref{5th-order1} the optimal linear weights are $c_0=\frac{2}{10},c_1=\frac{5}{10},c_2=\frac{3}{10}$.

If we replace the optimal linear weights $\{c_0, c_1, c_2\}$ with nonlinear weights $\{\omega_0, \omega_1, \omega_2\}$, a fifth-order finite volume compact WENO scheme can be obtained
\cite{guo2014positivity}
\begin{align}
\frac{2\omega_0+\omega_1}{3}u_{j-\frac{1}{2}}^{-}+&
\frac{\omega_0+2(\omega_1+\omega_2)}{3}u_{j+\frac{1}{2}}^{-}
+\frac{1}{3}\omega_2u_{j+\frac{3}{2}}^{-} \nonumber \\
=&\frac{1}{6}\omega_0\bar{u}_{j-1}+\frac{5(\omega_0+\omega_1)+\omega_2}{6}\bar{u}_j
+\frac{\omega_1+5\omega_2}{6}\bar{u}_{j+1}.
\label{CWENO-FV}
\end{align}
A set of nonlinear weights $\omega_k$ can be taken as \cite{castro2011high}
\begin{equation}
\label{weight-z}
\omega_k=\frac{\alpha_k^z}{\sum_{l=0}^2\alpha_l^z},\quad
\alpha_k^z=c_k\left(1+\left(\frac{\tau_5}{\beta^z_k+\epsilon}\right)^p
\right), \quad k=0,1,2,
\end{equation}
where $p \geq 1$ is the power parameter. $\epsilon$ is a small
positive number to avoid the denominator to be $0$. In our numerical
tests, we take $p=2$ and $\epsilon=10^{-13}$.

The smooth indicators $\beta_k^z$ are chosen from the WENO-Z scheme \cite{borges2008improved}
\begin{equation}
\beta_k^z=\left(\frac{\beta_k+\epsilon}{\beta_k+\tau_5+\epsilon}\right),\quad k=0,1,2,
\end{equation}
where $\tau_5=|\beta_{2}-\beta_{0}|$, which can improve the order of accuracy around the
smooth extrema as compared to the classic WENO scheme \cite{jiang1996efficient}.
The classical smooth indicators $\beta_k$ $(k=0,1,2)$ in \cite{jiang1996efficient} are given by
\begin{equation*}
\begin{aligned}
&\beta_{0}=\frac{13}{12}(\bar{u}_{j-2}-2\bar{u}_{j-1}+\bar{u}_{j})^2+\frac{1}{4}(\bar{u}_{j-2}-4\bar{u}_{j-1}+3\bar{u}_{j})^2,&\\
&\beta_{1}=\frac{13}{12}(\bar{u}_{j-1}-2\bar{u}_{j}+\bar{u}_{j+1})^2+\frac{1}{4}(\bar{u}_{j-1}-\bar{u}_{j+1})^2,&\\
&\beta_{2}=\frac{13}{12}(\bar{u}_{j}-2\bar{u}_{j+1}+\bar{u}_{j+2})^2+\frac{1}{4}(3\bar{u}_{j}-4\bar{u}_{j+1}+\bar{u}_{j+2})^2.&
\end{aligned}
\end{equation*}

We need to solve a tridiagonal system of \eqref{CWENO-FV} to get
$u_{j+\frac{1}{2}}^{-}$ at each time step since the nonlinear weights depend on the solutions.
$u_{j+\frac12}^{+}$ can be obtained similarly. See \cite{guo2014positivity} for
more discussions.

For the high-order compact scheme \eqref{CWENO-FV}, we need to set
appropriate boundary closures due to the global
nature of the reconstruction where all the flux values are involved
at each time step. Here for the scheme near boundaries, a fifth-order
WENO approximation is used \cite{ghosh2012compact}.

\subsection{Maximum-principle-satisfying limiter}
In \cite{guo2014positivity} a fifth-order finite volume compact WENO scheme
with positivity-preserving limiter was proposed for solving compressible Euler
equations. For scalar hyperbolic conservation laws \eqref{adveqn}, a similar idea with
the polynomial rescaling limiter \cite{zhang2010maximum,zhang2011maximum} for preserving the maximum principle is incorporated into the finite volume compact WENO scheme. The Euler forward temporal discretization for the semi-discrete scheme \eqref{1DadvFV-b} is
\begin{equation}\label{1DadvFV-c}
\bar{u}_j^{n+1}=\bar{u}_j^{n}-\lambda[h(u_{j+\frac{1}{2}}^{-},u_{j+\frac{1}{2}}^{+})
-h(u_{j-\frac{1}{2}}^{-},u_{j-\frac{1}{2}}^{+}) ],
\end{equation}
where $\lambda=\Delta t/\Delta x$. $u_{j+\frac{1}{2}}^{-}$ and
$u_{j+\frac{1}{2}}^{+}$ are the high order approximations of $u(x_{j+\frac{1}{2}},t^n)$
from the left and right limits which are reconstructed from
the finite volume compact WENO scheme. For simplicity, let $m=\min_x u(x,0)$
and $M=\max_x u(x,0)$, the polynomial rescaling limiter
proposed in \cite{zhang2011maximum} can be written as
\begin{equation}\label{MLimiterPPM}
\tilde{p}_j(x)=\theta(p_j(x)-\bar{u}_j)+\bar{u}_j, \quad
\theta=\min\Big\{\Big|\frac{M-\bar{u}_j^n}{M_j-\bar{u}_j^n}\Big|,\Big|\frac{m-\bar{u}_j^n}{m_j-\bar{u}_j^n}\Big|,1\Big\},
\end{equation}
with
\begin{equation}\label{Mm}
M_j=\max \{p_j(x_j^{*}),u_{j+\frac{1}{2}}^{-},
u_{j+\frac{1}{2}}^{+}\}, \quad m_j=\min \{p_j(x_j^{*}),u_{j+\frac{1}{2}}^{-},
u_{j+\frac{1}{2}}^{+}\}.
\end{equation}
$p_j(x)$ can be seen as a reconstructed polynomial with degree $4$ from
the cell-averaged values $\{\bar u_{j-1}, \bar u_{j}, \bar u_{j+1}\}$ and two
boundary values $\{u_{j-\frac{1}{2}}^{+}, u_{j+\frac{1}{2}}^{-}\}$ in cell $I_j$
for a fifth order approximation. There exists a point $x_j^*$ in cell $I_j$ such that
\begin{equation}
p_j(x_j^{*})=\frac{\bar{u}_j^n-\hat{\omega}_1u_{j-\frac{1}{2}}^{+}-
\hat{\omega}_G u_{j+\frac{1}{2}}^{-}}{1-2\hat{\omega}_1}.
\end{equation}
$\hat{\omega}_1$ and $\hat{\omega}_G$ are the first and last
weights of an $G$-point Gauss-Lobatto quadrature rule. Let
$\tilde{u}_{j-\frac{1}{2}}^{+}=\tilde{p}_j(x_{j-\frac{1}{2}})$ and
$\tilde{u}_{j+\frac{1}{2}}^{-}=\tilde{p}_j(x_{j+\frac{1}{2}})$,
we get a revised scheme of \eqref{1DadvFV-c}
\begin{equation}\label{1DadvFVPPM}
\bar{u}_j^{n+1}=\bar{u}_j^{n}-\lambda[h(\tilde{u}_{j+\frac{1}{2}}^{-},\tilde{u}_{j+\frac{1}{2}}^{+})
-h(\tilde{u}_{j-\frac{1}{2}}^{-},\tilde{u}_{j-\frac{1}{2}}^{+}) ].
\end{equation}
The scheme \eqref{1DadvFVPPM} satisfies a strict maximum principle
for scalar conservation laws under the CFL condition
\begin{equation}\label{cfl1d}
\lambda \alpha \leq \hat{\omega}_1.
\end{equation}
with a global Lax-Friedrichs flux $h(u,v)=\frac12[f(u)+f(v)-\alpha (v-u)]$,
where $\alpha=\max_u |f'(u)|$.
In the present compact
scheme, although $u_{j+\frac{1}{2}}^{-}$ and $u_{j-\frac{1}{2}}^{+}$
are obtained globally which are different from those in
\cite{zhang2010maximum}, the constructed polynomial $p_j(x)$ can be
seen locally. Thus this limiter does not destroy the high
order accuracy, see the proof in \cite{zhang2011maximum,zhang2010maximum}
for more details and see \cite{guo2014positivity} for more discussions.

\subsection{Temporal discretization}
In the present work, a strong stability preserving (SSP)
high order Runge-Kutta time discretizations \cite{gottlieb2009high} can be used to improve the
temporal accuracy in (\ref{1DadvFVPPM}). A third-order SSP
Runge-Kutta method is given as
\begin{equation}\label{RK3}
\begin{aligned}
u^{(1)}&=u^n+\Delta tL(u^n),&\\
u^{(2)}&=\frac{3}{4}u^n+\frac{1}{4}u^{(1)}+\frac{1}{4}\Delta tL(u^{(1)}),&\\
u^{n+1}&=\frac{1}{3}u^n+\frac{2}{3}u^{(2)}+\frac{2}{3}\Delta
tL(u^{(2)}),&
\end{aligned}
\end{equation}
where $L(u)$ is the spatial operator. For
a multi-stage SSP Runge-Kutta time method, the MPP limiter will be applied and
the tridiagonal system \eqref{CWENO-FV} will be solved at each stage
of each time step.

\section{Two dimensional case}
\label{sec3}
\setcounter{equation}{0}
\setcounter{figure}{0}
\setcounter{table}{0}

In this section, we consider the finite volume
compact WENO scheme for solving the two
dimensional conservation law
\begin{equation}\label{2Dadveqn}
u_t+f(u)_x+g(u)_y=0,
\end{equation}
on the domain $[a,b]\times[c,d]$ with rectangular meshes
\begin{equation}
a=x_{\frac{1}{2}}<x_{\frac{3}{2}}<\cdots<x_{N_x+\frac{1}{2}}=b, \quad
c=y_{\frac{1}{2}}<y_{\frac{3}{2}}<\cdots<y_{N_y+\frac{1}{2}}=d.
\end{equation}
Denoting $\Delta x=(b-a)/N_x$ and $\Delta y=(d-c)/N_y$ for uniform sizes,
the finite volume scheme for \eqref{2Dadveqn} on cell
$I_{ij}=[x_{i-\frac{1}{2}},x_{i+\frac{1}{2}}]\times
[y_{j-\frac{1}{2}},y_{j+\frac{1}{2}}]$ can be obtained as follows
\begin{equation}\label{2DadvFV-a}
\begin{aligned}
\bar{u}_{ij}^{n+1}=&\bar{u}_{ij}^{n} -\frac{\Delta t}{\Delta x
\Delta y}\int_{y_{j-\frac{1}{2}}}^{y_{j+\frac{1}{2}}}
[\hat{f}(u_{i+\frac{1}{2},j}^{-}(y),u_{i+\frac{1}{2},j}^{+}(y))-
\hat{f}(u_{i-\frac{1}{2},j}^{-}(y),u_{i-\frac{1}{2},j}^{+}(y))]dy\\
&-\frac{\Delta t}{\Delta x \Delta
y}\int_{x_{i-\frac{1}{2}}}^{x_{i+\frac{1}{2}}}
[\hat{g}(u_{i,j+\frac{1}{2}}^{-}(x),u_{i,j+\frac{1}{2}}^{+}(x))-
\hat{g}(u_{i,j-\frac{1}{2}}^{-}(x),u_{i,j-\frac{1}{2}}^{+}(x))]dx,
\end{aligned}
\end{equation}
where $u_{i+\frac{1}{2},j}^{-}(y)$, $u_{i-\frac{1}{2},j}^{+}(y)$,
 $u_{i,j+\frac{1}{2}}^{-}(x)$ and $u_{i,j-\frac{1}{2}}^{+}(x)$ denote
the traces of fifth-order polynomial on the four edges of cell
$I_{i,j}$ respectively. The cell-averaged value of cell $I_{i,j}$,
which denotes to be $\bar{u}_{ij}^{n}$, can be defined as
\begin{equation}
\bar{u}_{ij}(t)=\frac{1}{\Delta x \Delta y}\int_{x_{i-\frac{1}{2}}}^{x_{i+\frac{1}{2}}}
\int_{y_{j-\frac{1}{2}}}^{y_{j+\frac{1}{2}}}u(x,y,t)dxdy.
\end{equation}

We will use the Lax-Friedrichs flux for $\hat{f}$ and $\hat{g}$
\begin{subequations}
\label{LF}
\begin{align}
\hat{f}(u,v)=&\frac{1}{2}[f(u)+f(v)-a_1(v-u)],\quad a_1=\max_u|f'(u)|,\label{LFf}\\
\hat{g}(u,v)=&\frac{1}{2}[g(u)+g(v)-a_2(v-u)],\quad a_2=\max_u|g'(u)|.\label{LFg}
\end{align}
\end{subequations}
The integrals in \eqref{2DadvFV-a} are approximated by the Gaussian quadrature rule
as in \cite{zhang2010maximum} with sufficient accuracy. By using an $L$-point Gaussian quadrature
rule, we can get an approximation for \eqref{2DadvFV-a}
\begin{equation}\label{2DadvFV-b}
\begin{aligned}
\bar{u}_{ij}^{n+1}=&\bar{u}_{ij}^{n}
-\lambda_1\sum_{\beta=1}^L\omega_{\beta}
[\hat{f}(u_{i+\frac{1}{2},\beta}^{-},u_{i+\frac{1}{2},\beta}^{+})-
 \hat{f}(u_{i-\frac{1}{2},\beta}^{-},u_{i-\frac{1}{2},\beta}^{+})]\\
&-\lambda_2\sum_{\beta=1}^L\omega_{\beta}
[\hat{g}(u_{\beta,j+\frac{1}{2}}^{-},u_{\beta,j+\frac{1}{2}}^{+})-
 \hat{g}(u_{\beta,j-\frac{1}{2}}^{-},u_{\beta,j-\frac{1}{2}}^{+})],
 \end{aligned}
\end{equation}
where $\lambda_1=\frac{\Delta t}{\Delta x}$, $\lambda_2=\frac{\Delta
t}{\Delta y}$,
$u_{i\pm\frac{1}{2},\beta}^{\mp}=u_{i\pm\frac{1}{2},j}^{\mp}(y_j^{\beta})$,
$u_{\beta,j\pm\frac{1}{2}}^{\mp}=u_{i,j\pm\frac{1}{2}}^{\mp}(x_i^{\beta})$,
$y_j^{\beta}$ denotes the Gaussian quadrature points on
$[y_{j-\frac{1}{2}},y_{j+\frac{1}{2}}]$ and $x_i^{\beta}$ denotes
the Gaussian quadrature points on
$[x_{i-\frac{1}{2}},x_{i+\frac{1}{2}}]$, $\omega_{\beta}$ is
the corresponding quadrature weight.

The algorithm for maximum-principle-satisfying fifth-order finite volume
compact WENO scheme solving \eqref{2Dadveqn} consists of
the following steps:
\begin{enumerate}
\item The fifth-order finite volume
compact WENO scheme \eqref{CWENO-FV} is used to get
the four edge averages $\{\bar{u}_{i+\frac{1}{2},j}^{-},
\bar{u}_{i+\frac{1}{2},j}^{+}\}$ for fixed $j$ and
$\{\bar{u}_{i,j+\frac{1}{2}}^{-}, \bar{u}_{i,j+\frac{1}{2}}^{+}\}$
for fixed $i$, e.g.,
$\bar{u}_{i+\frac{1}{2},j}^{-}=\frac{1}{\Delta y}\int_{y_{j-\frac12}}^{y_{j+\frac12}}u(x_{i+\frac12}^-,y)dy$, similarly for others. Values at the quadrature points
$\{u_{i+\frac{1}{2},\beta}^{-}, u_{i+\frac{1}{2},\beta}^{-}\}$ and
$\{u_{\beta,j+\frac{1}{2}}^{-}, u_{\beta,j+\frac{1}{2}}^{-}\}$ are
obtained by the fifth-order WENO schemes \cite{jiang1996efficient}.
\item
Let $m=\min_{x,y} u(x,y,0)$ and $M=\max_{x,y} u(x,y,0)$, the maximum-principle-satisfying limiter
is constructed as follows \cite{zhang2011maximum}
\begin{equation}
\theta_{ij}=\min
\left\{\left|\frac{M-\bar{u}_{ij}^n}{M_{ij}-\bar{u}_{ij}^n}\right|
,\left|\frac{m-\bar{u}_{ij}^n}{m_{ij}-\bar{u}_{ij}^n}\right|,1\right\},
\end{equation}
where
\begin{equation}
M_{ij}=\max\{p_{ij}(x_i^{*},y_j^{*}),
u_{i\mp\frac{1}{2},\beta}^{\pm},u_{\beta,j\mp\frac{1}{2}}^{\pm}\}, \quad
m_{ij}=\min\{p_{ij}(x_i^{*},y_j^{*}),
u_{i\mp\frac{1}{2},\beta}^{\pm},u_{\beta,j\mp\frac{1}{2}}^{\pm}\},
\end{equation}
and there exists a point $(x_i^*, y_j^*)$ in cell $I_{ij}$ such that
\begin{equation}
p_{ij}(x_i^{*},y_j^{*})=
\frac{\bar{u}_{ij}^n-\sum_{\beta=1}^L\omega_{\beta}\hat{\omega}_1
[\mu_1(u_{i+\frac{1}{2},\beta}^{-}+u_{i-\frac{1}{2},\beta}^{+})
+\mu_2(u_{\beta,j+\frac{1}{2}}^{-}+u_{\beta,j-\frac{1}{2}}^{+})]}{1-2\hat{\omega}_1}.
\end{equation}
here $\mu_1=\lambda_1 a_1/(\lambda_1 a_1+\lambda_2 a_2)$ and $\mu_2=\lambda_2 a_2/(\lambda_1 a_1+\lambda_2 a_2)$ with $a_1$, $a_2$ defined in \eqref{LF}.
Note that $\omega_{\beta}$ is the weight of an $L$-point Gaussian quadrature rule
in \eqref{2DadvFV-b} and $\hat \omega_1$ is the first weight of an $G$-point
Gauss-Lobatto quadrature with $G=4$ for a fifth order scheme.
\item Finally $u_{i\mp\frac{1}{2},\beta}^{\pm}$ and $u_{\beta,j\mp\frac{1}{2}}^{\pm}$
in (\ref{2DadvFV-b}) are updated by
\begin{equation}\label{npv}
u_{i\mp\frac{1}{2},\beta}^{\pm}=\theta_{ij}
(u_{i\mp\frac{1}{2},\beta}^{\pm,old}-\bar{u}_{ij}^n)+\bar{u}_{ij}^n, \quad
u_{\beta,j\mp\frac{1}{2}}^{\pm}=\theta_{ij}
(u_{\beta,j\mp\frac{1}{2}}^{\pm,old}-\bar{u}_{ij}^n)+\bar{u}_{ij}^n.
\end{equation}
where $u_{i\mp\frac{1}{2},\beta}^{\pm,old}$ and $u_{\beta,j\mp\frac{1}{2}}^{\pm,old}$
are the values in (\ref{2DadvFV-b}) before modified.
\end{enumerate}
The CFL condition for the two dimensional case with the MPP limiter is taken to be
\begin{equation}\label{cfl2d}
a_1\lambda_1+a_2\lambda_2\leq \hat{\omega}_1.
\end{equation}
Similarly as the one dimensional case, the limiter does not destroy the high order accuracy.
The proof can be referred to \cite{zhang2010maximum}.

\section{Application to two dimensional incompressible flows}
\setcounter{equation}{0} \setcounter{figure}{0}
\setcounter{table}{0}
We now consider the incompressible Euler
equations in the vorticity stream-function formulation
\cite{bell1989second},
\begin{subequations}
\label{incom}
\begin{align}
&\omega_t+(u\omega)_x+(v\omega)_y=0, \label{incom-Euler}\\
&\Delta \psi =\omega, (u,v)=(-\phi_y,\phi_x), \label{incom-vor}
\end{align}
\end{subequations}
with $\omega(x,y,0)=\omega_0(x,y)$, $(u,v)\cdot\mathbf{n}=$ given on
$\partial \Omega$.

The divergence-free condition $u_x+v_y=0$ can be obtained from
equation \eqref{incom-vor}, which implies \eqref{incom-Euler} is
equivalent to the non-conservative form
\begin{equation}\label{incom-nonc}
\omega_t+u\omega_x+v\omega_y=0.
\end{equation}
The conservative equation \eqref{incom-Euler} itself does not imply the
maximum principle $\omega(x,y,t)\in [m, M]$ if without the incompressibility
condition $u_x+v_y=0$, where $m=\min_{x,y}\omega(x,y,0)$ and $M=\max_{x,y}\omega(x,y,0)$.
This is the main difficulty to get a maximum-principle-satisfying scheme for
solving \eqref{incom} \cite{zhang2010maximum}.
It has been proved in \cite{zhang2010maximum} that a high-order DG scheme
\cite{liu2000high} for \eqref{incom} with the MPP limiter under a suitable CFL
condition satisfies the maximum principle without destroying
the high order accuracy. In the following, we will consider the finite volume
compact WENO scheme on rectangular meshes similarly as in Section \ref{sec3}
for solving \eqref{incom-Euler}.

The finite volume scheme with Euler forward time discretization for
\eqref{incom-Euler} on cell $I_{i,j}$ is
\begin{equation}\label{incom-Euler-f}
\begin{aligned}
\bar{u}_{ij}^{n+1}=&\bar{u}_{ij}^{n} -\frac{\Delta t}{\Delta x
\Delta y}\int_{y_{j-\frac{1}{2}}}^{y_{j+\frac{1}{2}}}
[\hat{h}(\omega_{i+\frac{1}{2},j}^{-}(y),\omega_{i+\frac{1}{2},j}^{+}(y),u_{i+\frac{1}{2},j}(y))-
 \hat{h}(\omega_{i-\frac{1}{2},j}^{-}(y),\omega_{i-\frac{1}{2},j}^{+}(y),u_{i-\frac{1}{2},j}(y))]dy\\
&-\frac{\Delta t}{\Delta x \Delta
y}\int_{x_{i-\frac{1}{2}}}^{x_{i+\frac{1}{2}}}
[\hat{h}(\omega_{i,j+\frac{1}{2}}^{-}(x),\omega_{i,j+\frac{1}{2}}^{+}(x),v_{i,j+\frac{1}{2}}(x))-
 \hat{h}(\omega_{i,j-\frac{1}{2}}^{-}(x),\omega_{i,j-\frac{1}{2}}^{+}(x),v_{i,j-\frac{1}{2}}(x))]dx,
\end{aligned}
\end{equation}
where $\omega_{i+\frac{1}{2},j}^{-}(y)$,
$\omega_{i-\frac{1}{2},j}^{+}(y)$,
 $\omega_{i,j+\frac{1}{2}}^{-}(x)$ and $\omega_{i,j-\frac{1}{2}}^{+}(x)$ denote
the traces of a fifth-order polynomial on the four edges of cell
$I_{i,j}$ respectively.
%$u_{i+\frac{1}{2},j}^{-}(y)$,
%$u_{i-\frac{1}{2},j}^{+}(y)$,
% $v_{i,j+\frac{1}{2}}^{-}(x)$ and $v_{i,j-\frac{1}{2}}^{+}(x)$ denote the
% velocity on the left, right, bottom and top edges within cell $I_{i,j}$.
The cell-averaged value of the vorticity on cell $I_{i,j}$,
which denotes to be $\bar{\omega}_{ij}^{n}$, can be defined as
\begin{equation}
\bar{\omega}_{ij}(t)=\frac{1}{\Delta x \Delta
y}\int_{x_{i-\frac{1}{2}}}^{x_{i+\frac{1}{2}}}
\int_{y_{j-\frac{1}{2}}}^{y_{j+\frac{1}{2}}}\omega(x,y,t)dxdy.
\end{equation}
The Lax-Friedrichs upwind biased flux is used in this work, for instance on the right edge
\begin{equation}
\hat{h}(\omega_{i+\frac{1}{2},j}^{-}(y),\omega_{i+\frac{1}{2},j}^{+}(y),u_{i+\frac{1}{2},j}(y))
=\frac{1}{2}[u_{i+\frac{1}{2},j}(y)(\omega_{i+\frac{1}{2},j}^{-}(y)+\omega_{i+\frac{1}{2},j}^{+}(y))
-\alpha(\omega_{i+\frac{1}{2},j}^{+}(y)-\omega_{i+\frac{1}{2},j}^{-}(y))],
\end{equation}
where $\alpha$ is the maximum of $|u_{i+\frac{1}{2},j}(y)|$ either
locally or globally.

By using the $L$-point Gaussian quadrature rule, an approximation for
\eqref{incom-Euler-f} can be written as
\begin{equation}\label{incom-Euler-fvc}
\begin{aligned}
\bar{u}_{ij}^{n+1}=&\bar{u}_{ij}^{n}
-\lambda_1\sum_{\beta=1}^L\omega_{\beta}
[\hat{h}(\omega_{i+\frac{1}{2},\beta}^{-},\omega_{i+\frac{1}{2},\beta}^{+},
              u_{i+\frac{1}{2},\beta})-
 \hat{h}(\omega_{i-\frac{1}{2},\beta}^{-},\omega_{i-\frac{1}{2},\beta}^{+},
              u_{i-\frac{1}{2},\beta})]\\
&-\lambda_2\sum_{\beta=1}^L\omega_{\beta}
[\hat{h}(\omega_{\beta,j+\frac{1}{2}}^{-},\omega_{\beta,j+\frac{1}{2}}^{+},
              v_{\beta,j+\frac{1}{2}})-
 \hat{h}(\omega_{\beta,j-\frac{1}{2}}^{-},\omega_{\beta,j-\frac{1}{2}}^{+},
              v_{\beta,j-\frac{1}{2}})].
 \end{aligned}
\end{equation}
The algorithm for maximum-principle-satisfying fifth-order
finite volume compact WENO scheme in Section 3 can be
applied to \eqref{incom-Euler-fvc}. Under the CFL condition
$a_1\lambda_1+a_2\lambda_2\leq \frac{1}{2}\min_{\alpha=1,\cdots
,G}\hat{\omega}_{\alpha}$, the scheme
\eqref{incom-Euler-fvc} satisfies the maximum principle and maintains the high
order accuracy \cite{zhang2010maximum}. For
incompressible Euler equations in the vorticity stream-function
formulation \eqref{incom}, the Possion equation \eqref{incom-vor} is solved
by the Fourier spectral method.

\section{Numerical examples}
\setcounter{equation}{0}
\setcounter{figure}{0}
\setcounter{table}{0}

In this section, we provide some classical numerical examples for
the fifth-order finite volume compact WENO scheme with the third
order SSP Runge-Kutta time discreization \eqref{RK3}, which is denoted as the ``FVCW'' scheme.
The original fifth order WENO scheme of Jiang and Shu \cite{jiang1996efficient}
is denoted as ``WENO-JS'' and the fifth-order WENO-Z scheme
\cite{castro2011high} will be denoted as ``WENO-Z''. We will compare the
present FVCW scheme with WENO-JS and WENO-Z schemes in some of the following test cases.
We compute the solutions up to time $T$ on a mesh of $N$ and $N\times N$
with the CFL conditions (\ref{cfl1d}) and (\ref{cfl2d}) for one and
two dimensional cases respectively. The minimum and maximum numerical cell-averaged values are denoted as
$(\bar u_h)_{min}$ and $(\bar u_h)_{max}$ or $(\bar \omega_h)_{min}$ and $(\bar \omega_h)_{max}$ for the
incompressible flow problems respectively.

\subsection{One-dimensional test cases}
\begin{exa}\label{Eg:1}
We first solve the linear advection equation
\begin{equation}
u_t+u_x=0, \quad u(x,0)=u_0(x),
\label{ex1}
\end{equation}
with periodic boundary conditions on the domain $[0,2]$. We take the
smooth initial data $u_0(x)=0.5+\sin^4(2\pi x)$ to test the accuracy
and the maximum principle preserving property. In order to compare
our numerical results with those obtained by the finite volume WENO
schemes in \cite{zhang2011maximum}, the weights in
\cite{jiang1996efficient} are used. The $L^1$ and $L^{\infty}$
errors and orders at time $T=0.1$ with and
without limiters are given in Table \ref{tab:1}. We observe that
the numerical solutions obtained with limiters are all lied in
$[0.5,1.5]$, while the minimum values might be less than 0.5 if
without limiters. It shows the present FVCW scheme with limiters
satisfies the strict discrete maximum principle and the high order
of accuracy is maintained. The numerical results are comparable to
those obtained by the finite volume WENO scheme in Table 5.1 in
\cite{zhang2010maximum}.

\begin{table}
\centering
\caption{$L^1$ and $L^{\infty}$ errors and orders for
Example \ref{Eg:1} with $u_0(x)=0.5+\sin^4(2\pi x)$.}
\label{tab:1}
\begin{tabular}{lllllll}
\hline\noalign{\smallskip}
N  &  $L^1$ error & Order&$L^{\infty}$ error& Order & $(\bar u_h)_{min}$ &   $(\bar u_h)_{max}$ \\
\noalign{\smallskip} \hline \noalign{\smallskip}
\multicolumn{7}{c}{\textrm{with limiters}} \\
\noalign{\smallskip} \hline \noalign{\smallskip}
    20 &    5.56E-03 & &    1.23E-02 & &    0.5056844053    &    1.4336458051         \\
    40 &    9.02E-04 &     2.62 &    2.95E-03 &     2.06 &    0.5005823859   &     1.4838817402        \\
    80 &    5.28E-05 &     4.09 &    2.68E-04 &     3.46 &    0.5000324661   &     1.4959125035        \\
   160 &    8.81E-07 &     5.91 &    6.96E-06 &     5.27 &    0.5000001415   &     1.4989730706        \\
   320 &    2.10E-08 &     5.39 &    1.19E-07 &     5.87 &    0.5000000287   &     1.4997430328        \\
   640 &    4.82E-10 &     5.44 &    1.75E-09 &     6.10 &    0.5000000018   &     1.4999357475        \\
\noalign{\smallskip} \hline \noalign{\smallskip}
   \multicolumn{7}{c}{\textrm{without limiters}} \\
\noalign{\smallskip} \hline \noalign{\smallskip}
    20 &    4.80E-03 & &    1.38E-02 & &    0.4988356090    &    1.4364151979         \\
    40 &    1.01E-03 &     2.24 &    4.04E-03 &     1.77 &    0.4996159401    &    1.4844567924         \\
    80 &    6.73E-05 &     3.91 &    3.35E-04 &     3.59 &    0.4999199982    &    1.4959203562         \\
   160 &    8.83E-07 &     6.25 &    6.89E-06 &     5.60 &    0.4999999714    &    1.4989731984          \\
   320 &    2.11E-08 &     5.39 &    1.20E-07 &     5.85 &    0.5000000279    &    1.4997430383          \\
   640 &    4.82E-10 &     5.45 &    1.75E-09 &     6.10 &    0.5000000018    &    1.4999357479         \\
\noalign{\smallskip}\hline
\end{tabular}
\end{table}
\end{exa}

\begin{exa}\label{Eg:2}
We then consider the linear advection equation \eqref{ex1} with the following
initial condition \cite{jiang1996efficient}
\begin{equation}\label{ic:2}
u(x,0)=\left\{\begin{array}{ll}
\frac{1}{6}(G(x,z-\delta)+G(x,z+\delta)+4G(x,z)), \quad &
\textrm{ $x\in(-0.8,-0.6)$},\\
1, \quad & \textrm{ $x\in(-0.4,-0.2)$},\\
1-|10(x-0.1)|, \quad & \textrm{ $x\in(0,0.2)$},\\
\frac{1}{6}(F(x,a-\delta)+F(x,a+\delta)+4F(x,a)), \quad&
\textrm{ $x\in(0.4,0.6)$},\\
0, \quad &\textrm{otherwise},
\end{array}\right.
\end{equation}
where $G(x,z)=e^{-\beta(x-z)^2},
F(x,\gamma)=\sqrt{\max(1-\alpha^2(x-\gamma)^2,0)}$, $z=-0.7, 
\delta=0.0005, \alpha=10, \beta=\frac{\log 2}{36\delta^2},\gamma=0.5$. 
The computational
domain is $[-1,1]$ with periodic boundary conditions.

The solution is a travelling wave formed by the combination of a
Gaussian, a square wave, a sharp triangle wave and a half ellipse.
After a period of $T=2$, the solution will get back to its initial
position. In Fig. \ref{Fig:1}, we show the solution at time $T = 8$
on a grid with $N=200$ for the WENO-JS, WENO-Z and FVCW schemes. In
the zoom-in figures, we can see the FVCW scheme captures better
results that the other two schemes.

The maximum and minimum numerical solutions are listed in Table \ref{tab:2}. The
numerical solutions with limiters are all within the range $[0,1]$.
However, without limiters the minimum values are negative and the
maximum values are great than 1.

\begin{figure}
\subfigure[T=8] {\includegraphics[width=0.5\textwidth]{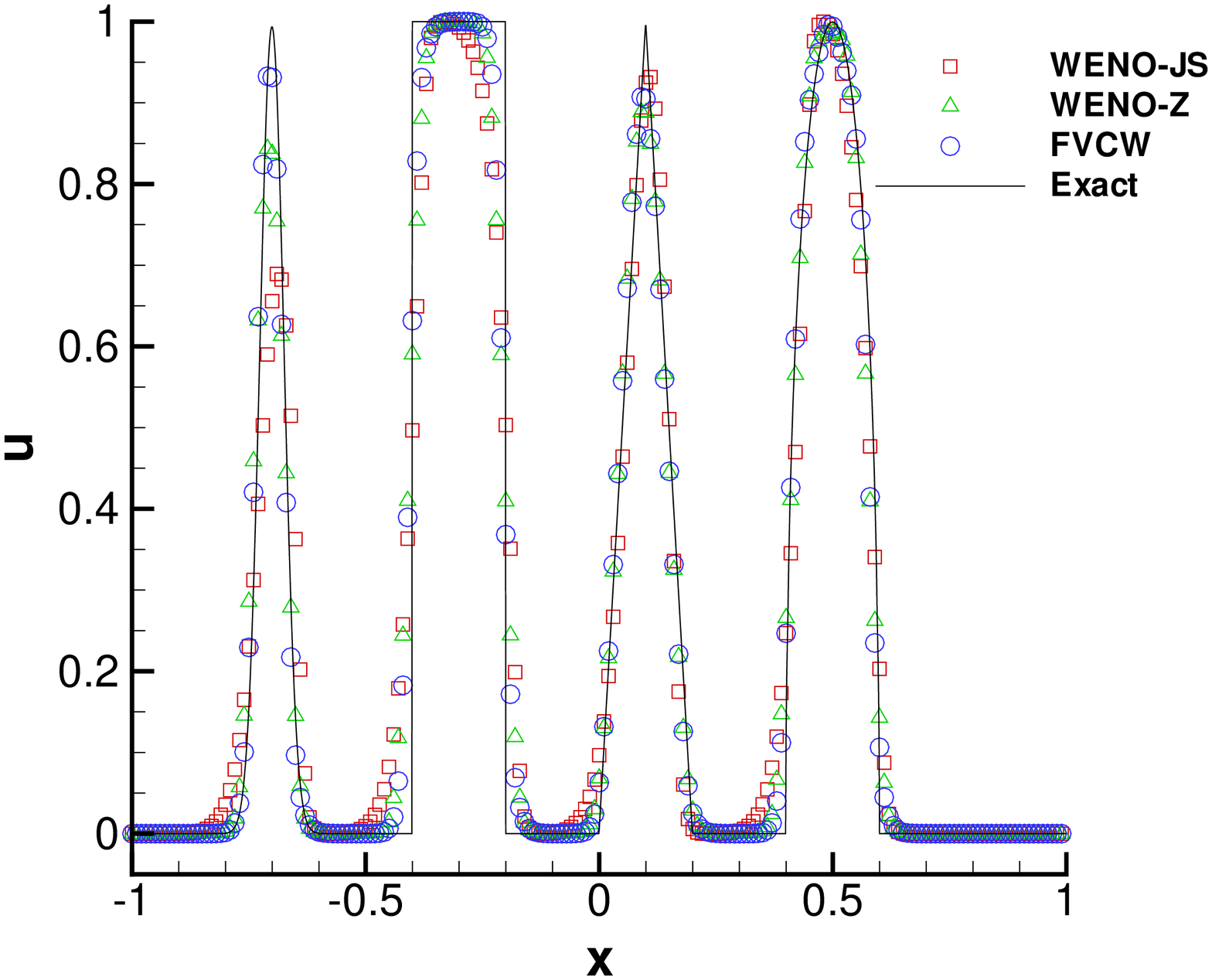}}
\subfigure[Gaussian
wave]{\includegraphics[width=0.5\textwidth]{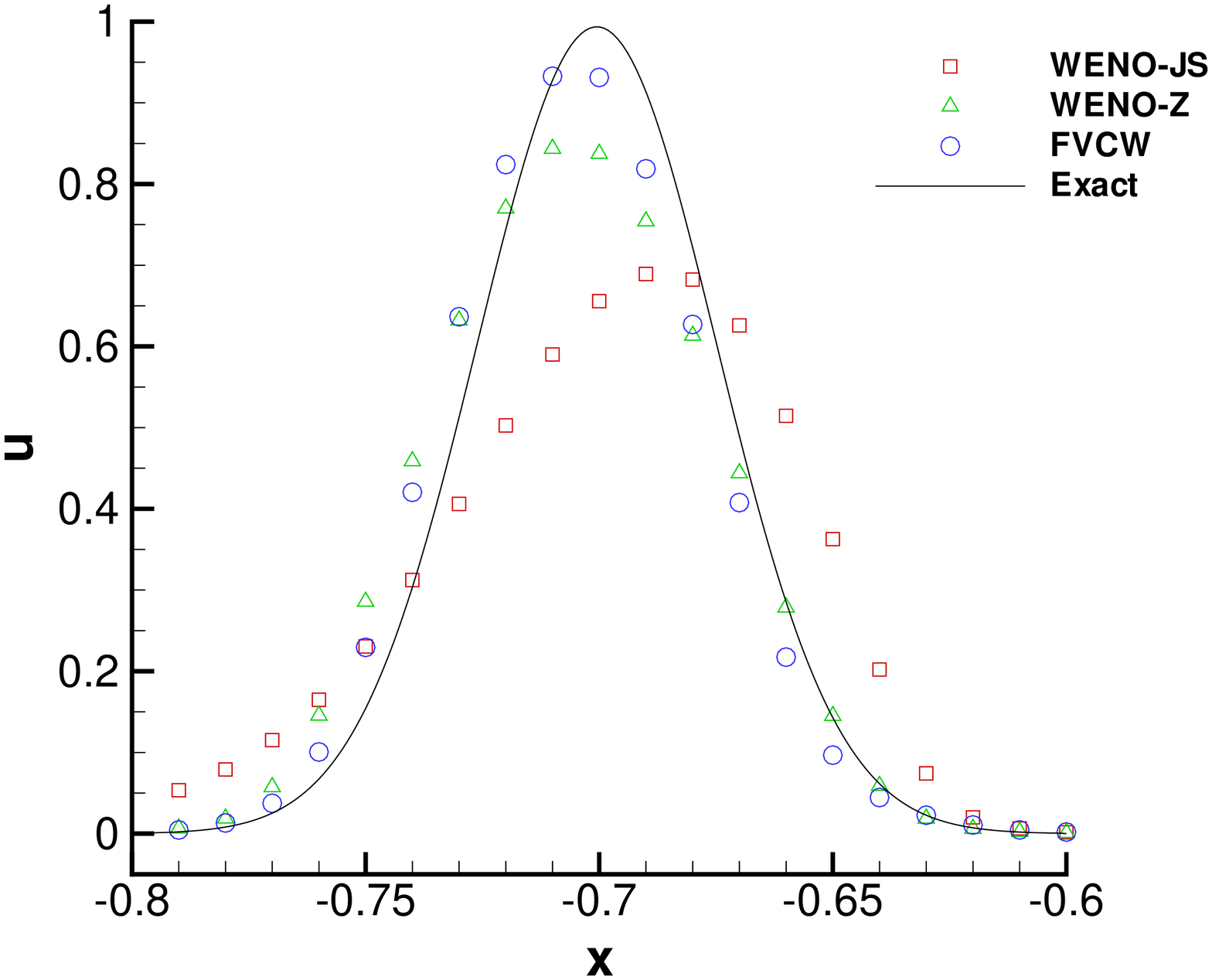}}
\subfigure[Square
wave]{\includegraphics[width=0.5\textwidth]{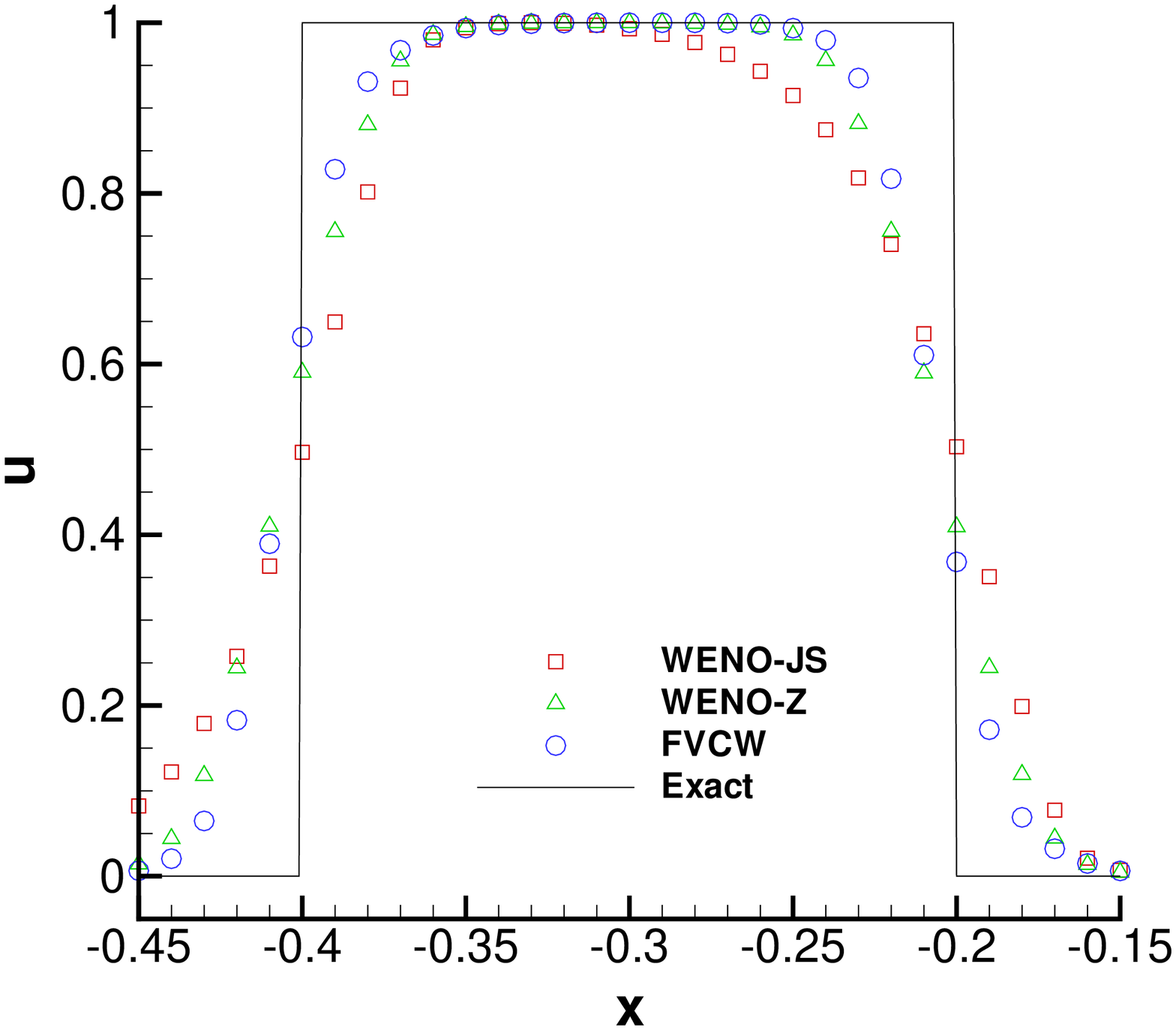}}
\subfigure[Triangular
wave]{\includegraphics[width=0.5\textwidth]{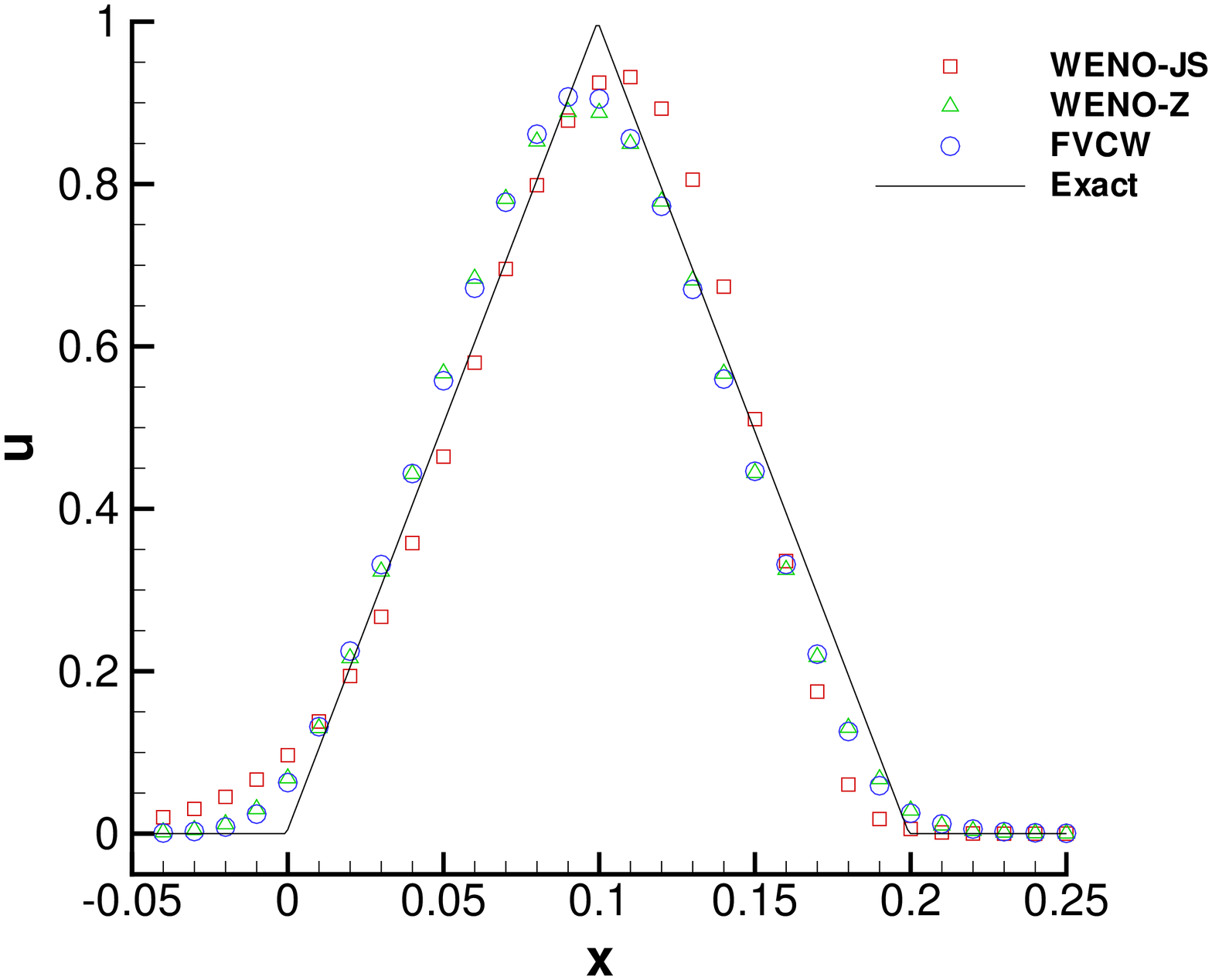}}
\caption{Numerical results computed by the WENO-JS, WENO-Z and FVCW
schemes with the exact solution for Example \ref{Eg:2} at $T=8$.
$N=200$.}\label{Fig:1}
\end{figure}

\begin{table}
\centering
\caption{Maximum and minimum numerical solutions for Example
\ref{Eg:2} at $T = 8$.} \label{tab:2}
\begin{tabular}{c|cc|cc}
\hline
   &  \multicolumn{2}{c|}{\textrm{with limiters}}  & \multicolumn{2}{c}{\textrm{without limiters}}\\
\hline
N  &  $(\bar u_h)_{min}$ & $(\bar u_h)_{max}$ &$(\bar u_h)_{min}$ & $(\bar u_h)_{max}$ \\
\hline
    50 &    6.2532042145E-05  &   0.9844087972 &-6.4754491922E-02   &  1.0136701866\\
   100 &    6.2444287122E-08  &   0.9985999188 &-4.1188418939E-03   &  1.0332676399\\
   200 &    5.9311237190E-13  &   0.9999770077 &-9.5377767652E-05   &  1.0043551468\\
   400 &    6.6787958057E-19  &   0.9999998937 &-1.2991894432E-07   &  1.0000004916\\
   800 &    6.0843658046E-31  &   1.0000000000 &-9.7496437457E-08   &  1.0000000859\\
   \hline
\end{tabular}
\end{table}
\end{exa}

\begin{exa}\label{Eg:3}
This example is the linear advection equation (\ref{ex1}) with initial condition \cite{shu1988efficient},
\begin{equation}\label{ic:3}
u(x+0.5,0)=\left\{\begin{array}{ll}
-x\sin(\frac{3}{2}\pi x^2), \quad & \textrm{$-1\leq x < \frac{1}{3}$},\\
|\sin(2\pi x)|, \quad & \textrm{$|x| \le \frac{1}{3}$},\\
2x-1-\frac{1}{6}\sin(3\pi x), \quad & \textrm{$\frac{1}{3}<x\le 1$},
\end{array}\right.
\end{equation}
on the domain $[-1,1]$ with periodic boundary conditions. This test case
is used to check how the FVCW scheme can capture smooth and discontinuous solution structures.
In Fig. \ref{Fig:2}, we show the solution at time $T = 2$ on the
grid with $N=200$ for the WENO-JS, WENO-Z and FVCW schemes. From the
figure we can see the numerical solutions all capture the
discontinuities without oscillations and  match the accurate solution
in the smooth regions very well. In the zoom-in Fig. \ref{Fig:2}(b), it
shows that the FVCW scheme performs better than the other
two schemes.

\begin{figure}
\subfigure[N=200]{\includegraphics[width=0.5\textwidth]{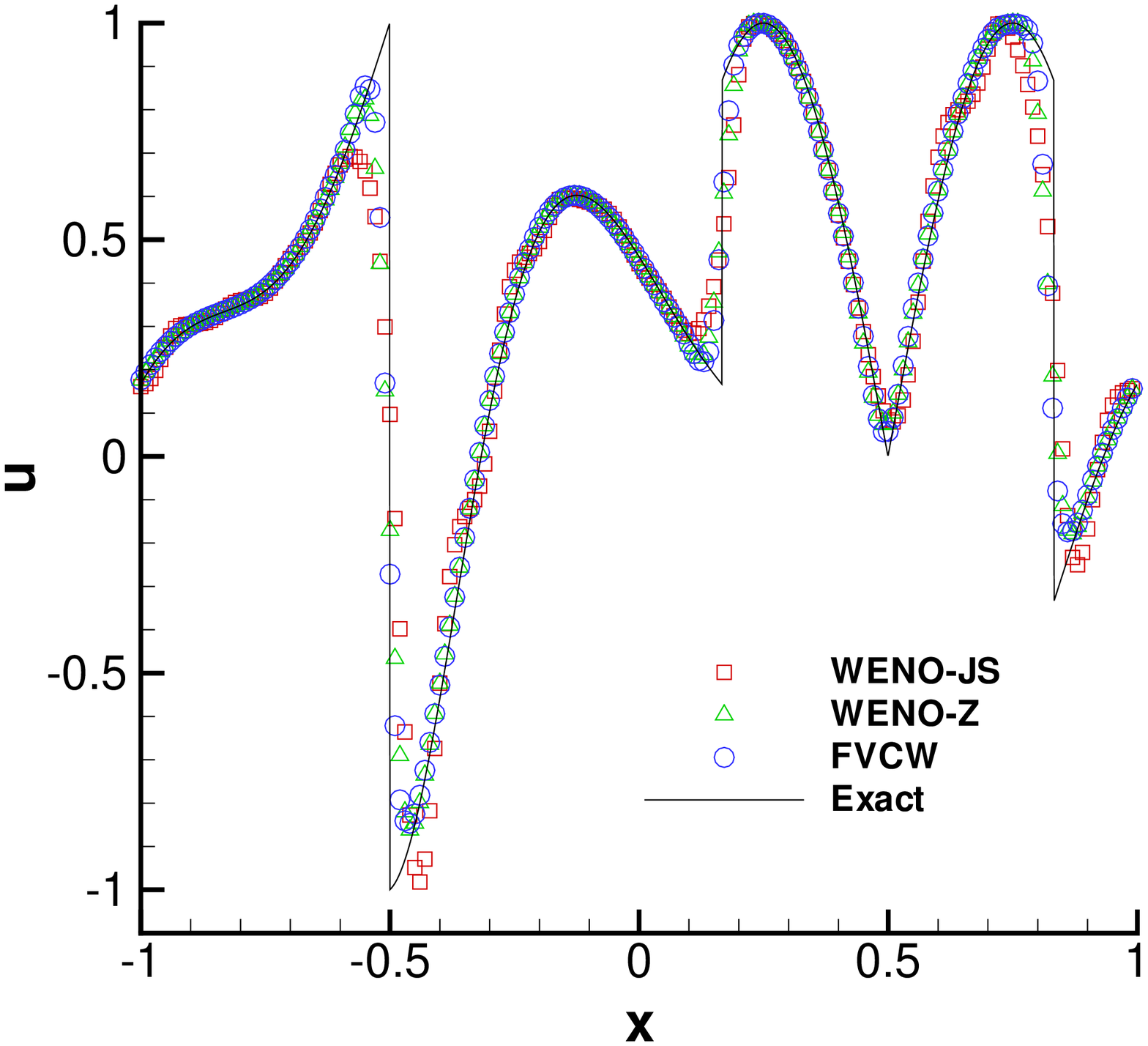}}
\subfigure[The enlarged portion of
(a).]{\includegraphics[width=0.5\textwidth]{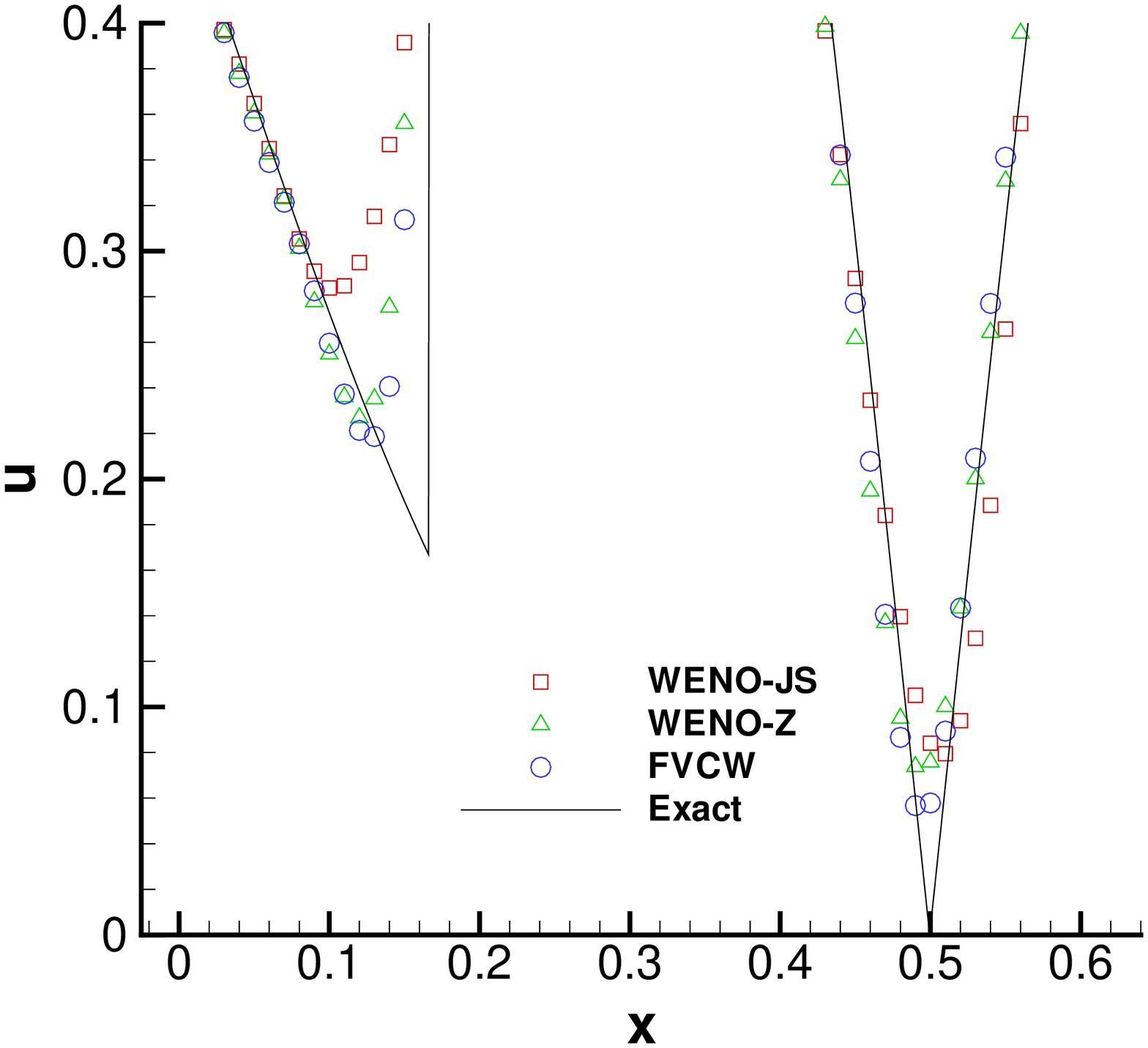}}
\caption{Numerical results computed by the WENO-JS, WENO-Z and FVCW
schemes with the exact solution for Example \ref{Eg:3} at $T=2$.
$N=200$.}\label{Fig:2}
\end{figure}
\end{exa}

\begin{exa}\label{Eg:4}
We now solve the nonlinear Burgers' equation
\begin{equation}
u_t+(\frac{u^2}{2})_x=0, \quad u(x,0)=u_0(x),
\end{equation}
with periodic boundary conditions. We test our scheme with the initial condition
$u_0(x)=\sin^{4}(x)$ on $[0,2\pi]$. At $T = 0.5$, the solution is smooth.
The $L^1$ and $L^{\infty}$ errors and orders
for the FVCW scheme with and without limiters are given in
Tables \ref{tab:3}. For this nonlinear problem, the numerical solutions with limiters are all
within the range $[-0.5,1.5]$ and the 5th order of accuracy is also maintained.
At $T = 1.2$ the solution develops a shock. The numerical solution
and the exact solution are shown in Fig. \ref{Fig:3}. With limiters, the minimum numerical value
is 5.25E-006. However, it is -2.84E-006 if without limiters.

\begin{table}
\centering
\caption{$L^1$ and $L^{\infty}$ errors and orders for
Example \ref{Eg:4} with $u_0(x)=\sin^4(x)$.}
\label{tab:3}
\begin{tabular}{lllllll}
\hline\noalign{\smallskip}
N  &  $L^1$ error & Order&$L^{\infty}$ error& Order & $(\bar u_h)_{min}$ & $(\bar u_h)_{max}$   \\
\noalign{\smallskip} \hline \noalign{\smallskip}
\multicolumn{7}{c}{\textrm{with limiters}} \\
\noalign{\smallskip} \hline \noalign{\smallskip}
    20 &     7.11E-03 &   &1.75E-02 &                      &4.8027078317E-03 &0.9661788333\\
    40 &     1.54E-03 &     2.21 &     9.32E-03 &     0.91 &5.9499022915E-04 &0.9860591424\\
    80 &     1.79E-04 &     3.11 &     2.08E-03 &     2.17 &5.3042057661E-05 &0.9984210631\\
   160 &     9.62E-06 &     4.21 &     1.88E-04 &     3.47 &4.2567524841E-06 &0.9995883857\\
   320 &     3.54E-07 &     4.76 &     7.63E-06 &     4.62 &3.2493209900E-07 &0.9999342274\\
   640 &     1.10E-08 &     5.00 &     2.36E-07 &     5.02 &1.6373928322E-08 &0.9999490507\\
\noalign{\smallskip} \hline \noalign{\smallskip}
\multicolumn{7}{c}{\textrm{without limiters}} \\
\noalign{\smallskip} \hline \noalign{\smallskip}
    20 &     6.34E-03 & &     1.85E-02 & &                   1.7218752333E-04 & 0.9664450834\\
    40 &     1.48E-03 &     2.10 &     9.06E-03 &     1.03 &-3.1070308046E-05 & 0.9860690299\\
    80 &     1.83E-04 &     3.02 &     2.08E-03 &     2.12 &-1.0538041762E-05 & 0.9984240105\\
   160 &     9.57E-06 &     4.25 &     1.88E-04 &     3.47 &-1.4739089684E-06 & 0.9995889053\\
   320 &     3.51E-07 &     4.77 &     7.63E-06 &     4.62 &-1.6598308446E-07 & 0.9999343533\\
   640 &     1.14E-08 &     4.95 &     2.36E-07 &     5.02 &-1.6619301474E-08 & 0.9999490606\\
\noalign{\smallskip}\hline
\end{tabular}
\end{table}

\begin{figure}[htp]
\subfigure[without
limiters]{\includegraphics[width=0.5\textwidth]{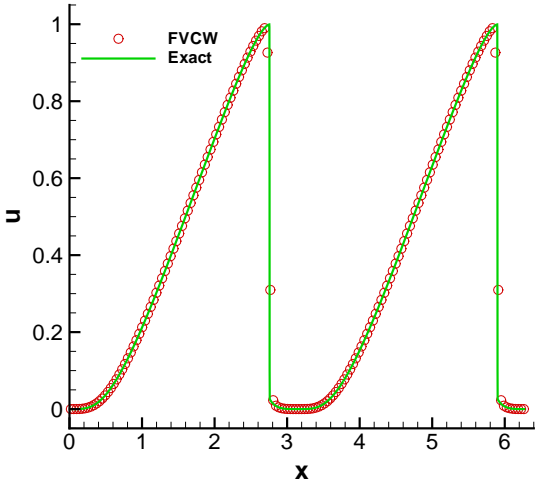}}
\subfigure[with
limiters]{\includegraphics[width=0.5\textwidth]{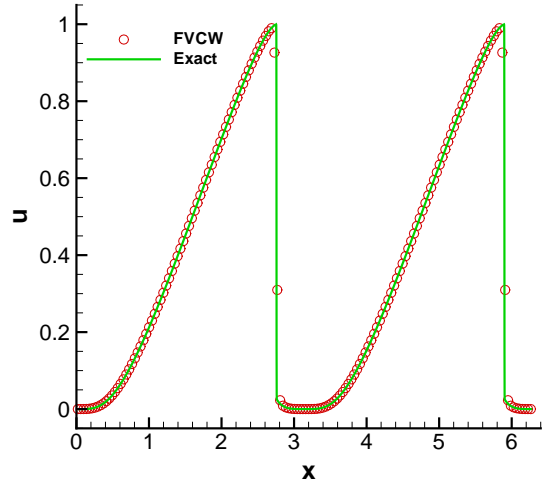}}
\caption{Numerical results with the exact solution for Example
\ref{Eg:4} at $T=1.2$ with $N=160$.}\label{Fig:3}
\end{figure}
\end{exa}

\begin{exa}\label{Eg:5}
The nonlinear Buckley-Leverett problem \cite{shu1988efficient}
is used for reservoir simulation,
\begin{equation}
u_t+f(u)_x=0, \quad f(u)=\frac{4u^2}{4u^2+(1-u)^2}.
\label{ex5}
\end{equation}
$f(u)$ is a nonconvex function and the initial condition is
\begin{equation}\label{ic:5}
u(x,0)=\left\{\begin{array}{ll}
1, \quad & \textrm{$-\frac{1}{2}< x <0$},\\
0, \quad & \textrm{otherwise}.
\end{array}\right.
\end{equation}

We compute the solution up to $T=0.4$ on the domain $[-1, 1]$ with
inflow and outflow boundary conditions on each side respectively.
The numerical solutions of WENO-JS, WENO-Z, FVCW with $N=100$ and
the exact solution are shown in Fig. \ref{Fig:4}. The zoom-in Fig.
\ref{Fig:4}(b) shows that a slightly more accurate result can be
obtained by the FVCW scheme. The maximum and minimum numerical
solutions are listed in Table \ref{tab:4} for the FVCW scheme with
and without limiters. The numerical solutions with limiters are all
within the range $[0,1]$ as expected.

\begin{figure}[htp]
\subfigure[N=100]{\includegraphics[width=0.5\textwidth]{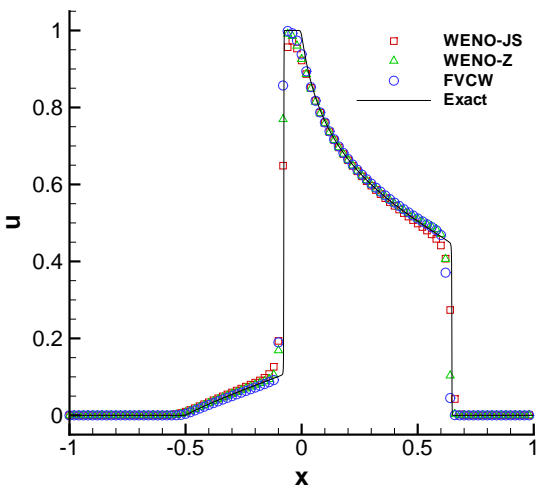}} %N=100
\subfigure[The enlarged portion of (a).]
{\includegraphics[width=0.5\textwidth]{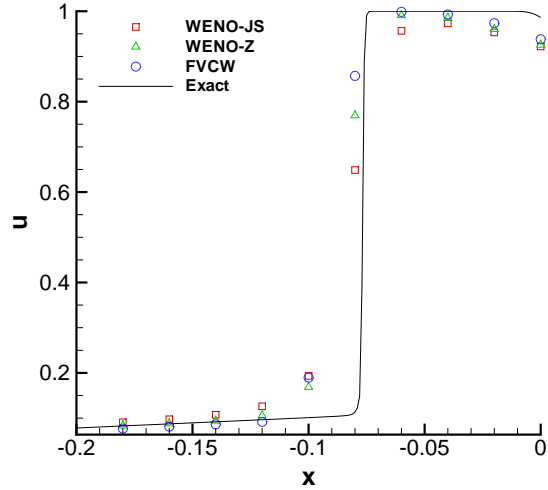}} %N=200
\caption{Numerical results computed by the WENO-JS, WENO-Z and FVCW
with the exact solution for Example \ref{Eg:5} at $T=0.4$ with
$N=100$.}\label{Fig:4}
\end{figure}

\begin{table}
\centering
\caption{Maximum and minimum numerical solutions for Example
\ref{Eg:5} at $T = 0.4$.} \label{tab:4}
\begin{tabular}{c|cc|cc}
\hline%\noalign{\smallskip}
&  with limiters    &               & without limiters\\
\hline%\noalign{\smallskip}
N  &  $(\bar u_h)_{min}$ & $(\bar u_h)_{max}$ &$(\bar u_h)_{min}$ & $(\bar u_h)_{max}$ \\
%\noalign{\smallskip}
\hline
    50 &    2.3624923015E-13  &   0.9000000000 &-1.0004055700E-08   &  0.9000000000\\
   100 &    5.3108644691E-21  &   0.9876234859 &-1.2483926241E-08   &  0.9876241852\\
   200 &    1.9526892745E-30  &   0.9991089796 &-1.5844261578E-08   &  0.9991089584\\
   400 &    1.3188489112E-47  &   0.9999947633 &-1.8249582807E-08   &  0.9999947651\\
   800 &    7.5182650962E-74  &   1.0000000000 &-1.7783953062E-08   &  1.0000000166\\
\hline
\end{tabular}
\end{table}
\end{exa}

\subsection{Two-dimensional test cases}
\begin{exa}\label{Eg:6}
We first solve the two dimensional linear advection equation
\begin{equation}
u_t+u_x+u_y=0, \quad u(x,y,0)=u_0(x,y), \quad (x,y)\in[0,1]\times[0,1],
\end{equation}
with periodic boundary condition and a smooth initial data
$u_0(x,y)=\sin^4(2\pi(x+y))$. We compute the solution up to time $t
=0.1$. The $L^1$ and $L^{\infty}$ errors and orders for the FVCW
scheme with and without limiters are given in Table \ref{tab:5}. The
FVCW scheme with limiters satisfies the strict maximum principle and
the 5th order of accuracy is maintained, which is similar to the one
dimensional case.

\begin{table}
\centering
\caption{$L^1$ and $L^{\infty}$ errors and orders for Example
\ref{Eg:6} with $u_0(x,y)=\sin^4(2\pi(x+y))$.} \label{tab:5}
\begin{tabular}{lllllllll}
\hline\noalign{\smallskip}
N $\times$ N  &  $L^1$ error & Order&$L^{\infty}$ error& Order & $(\bar u_h)_{min}$ & $(\bar u_h)_{max}$   \\
\noalign{\smallskip} \hline \noalign{\smallskip}
   \multicolumn{7}{c}{\textrm{with limiters}} \\
\noalign{\smallskip} \hline \noalign{\smallskip}
    10  $\times$  10 &      5.81E-02 & &     1.09E-01 & &3.2175761831E-02   &  0.9000000000\\
    20  $\times$  20 &     4.24E-03 &     3.78 &     1.33E-02 &     3.03
    &3.0806089368E-03  &   0.9611971380\\
    40  $\times$  40 &     5.78E-04 &     2.87 &     2.25E-03 &     2.57
    &3.2800792678E-04  &   0.9916153520\\
    80  $\times$  80 &     2.90E-05 &     4.32 &     2.21E-04 &     3.35
    &2.4330912515E-05  &   0.9979365718\\
   160  $\times$ 160 &     5.32E-07 &     5.77 &     3.78E-06 &     5.87
    &2.1151199593E-07  &   0.9994852045\\
    320  $\times$ 320 &     1.00E-08 &     5.73 &     1.01E-07 &     5.23
    &1.1731843775E-08   &  0.9998714151\\
\noalign{\smallskip} \hline \noalign{\smallskip}
   \multicolumn{7}{c}{\textrm{without limiters}} \\
\noalign{\smallskip} \hline \noalign{\smallskip}
    10  $\times$  10 &     5.17E-02 & &     8.22E-02 &          &-4.6206453729E-02    & 0.9000000000\\
    20  $\times$  20 &     4.10E-03 &     3.66 &     1.50E-02 &     2.46 &-1.0767640535E-02   &  0.9652658048\\
    40  $\times$  40 &     4.78E-04 &     3.10 &     1.72E-03 &     3.12 &-6.4056253913E-04   &  0.9918385835\\
    80  $\times$  80 &     2.84E-05 &     4.07 &     2.29E-04 &     2.91 &-7.1932689350E-05  &   0.9979464519\\
   160  $\times$ 160 &     5.05E-07 &     5.81 &     3.78E-06 &     5.92 &-3.2339904232E-08  &   0.9994861246\\
      320  $\times$ 320 &     8.04E-09 &     5.97 &     3.27E-08 &     6.85 &5.9884217874E-09  &   0.9998715025\\
\noalign{\smallskip}\hline
\end{tabular}
\end{table}

\end{exa}

\begin{figure}[htp]
\subfigure[without
limiters]{\includegraphics[width=0.5\textwidth]{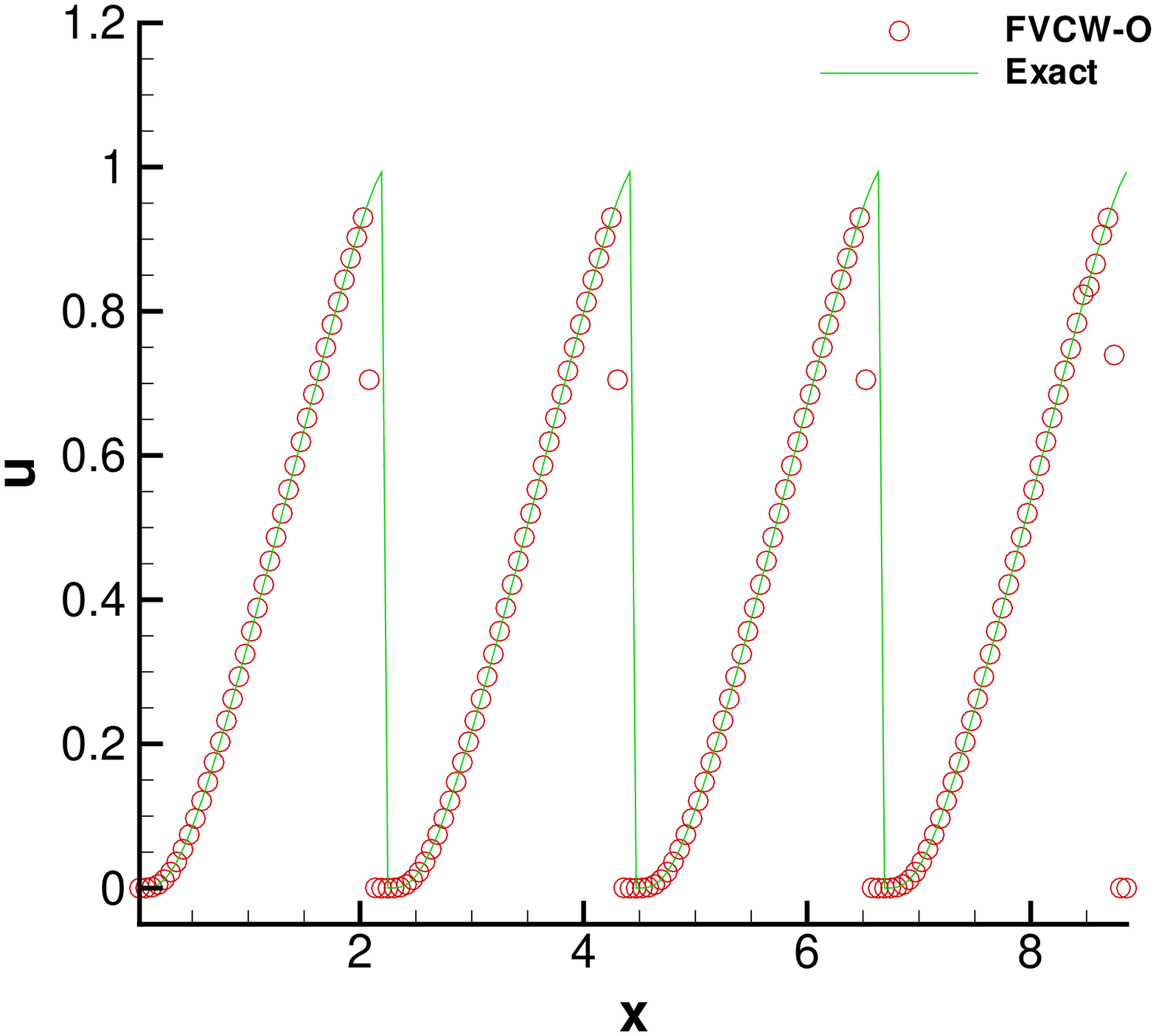}}
\subfigure[with
limiters]{\includegraphics[width=0.5\textwidth]{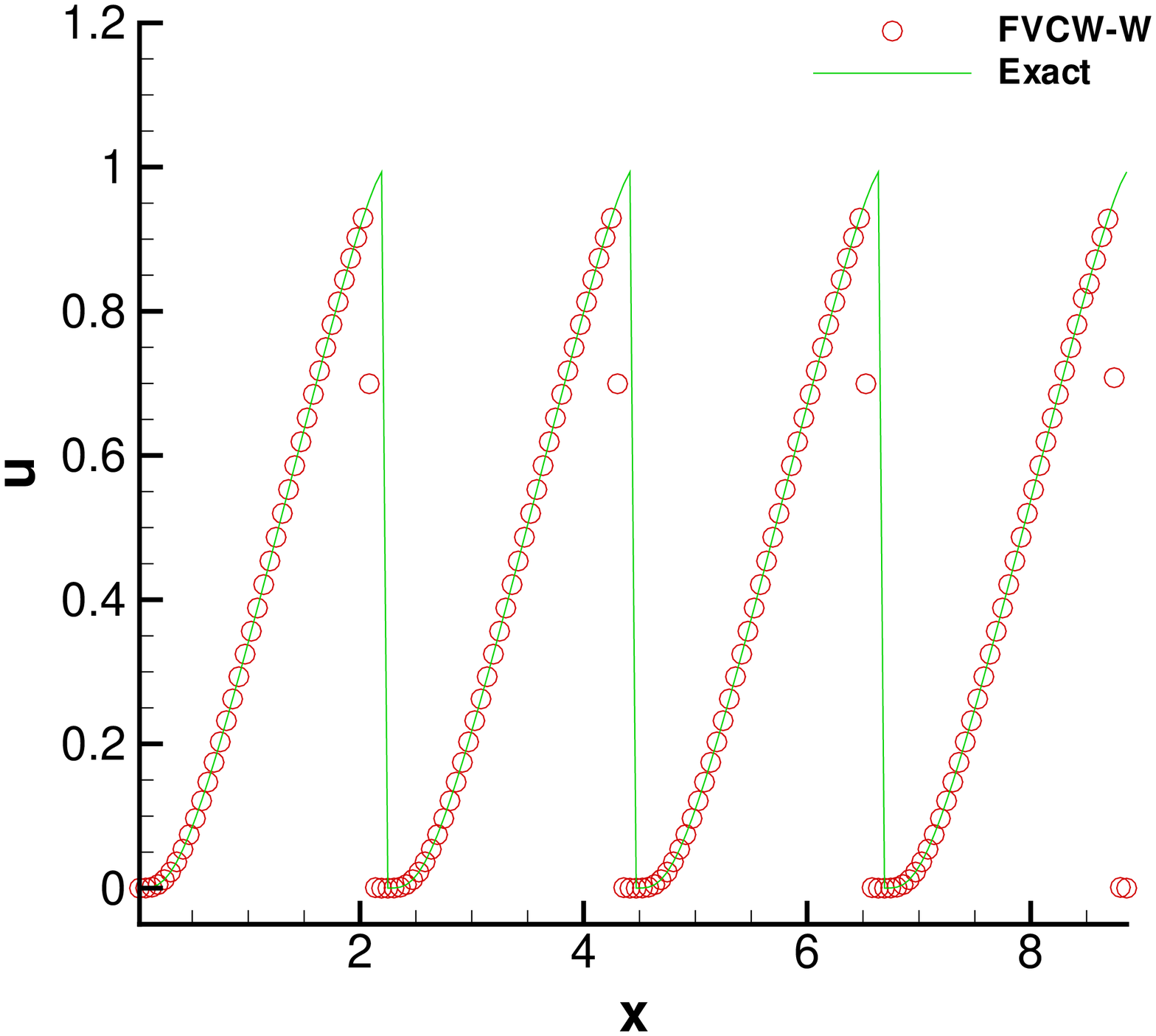}}
\caption{Numerical results for Example \ref{Eg:7} with
$u_0(x,y)=\sin^4(x+y)$ at $T=0.8$. Mesh $160 \times 160$. Cuts along
$x=y$.}\label{Fig:5}
\end{figure}

\begin{exa}\label{Eg:7}
We solve the nonlinear Burgers' equation
\begin{equation}
u_t+(\frac{u^2}{2})_x+(\frac{u^2}{2})_y=0,\quad u(x,y,0)=u_0(x,y),\quad (x,y)\in[0,2\pi]\times [0,2\pi],
\end{equation}
with periodic boundary condition.
\end{exa}
We also test our scheme with the initial condition
$u_0(x,y)=\sin^{4}(x+y)$ on $[0,2\pi]\times [0,2\pi]$ with periodic
boundary conditions.  We compute up to $T = 0.2$ when the solution
is still smooth.  The $L^1$ and $L^{\infty}$ errors and orders of
accuracy for the FVCW scheme with and without limiters are given in
Table \ref{tab:6}, almost fifth order accuracy can be observed for this example. At time $T =0.8$ when the
solution develops a still shock, the numerical solutions with $160
\times 160$ grid points and the exact solution are showed in Fig.
\ref{Fig:5}. With limiters, the minimum numerical value is 1.02E-006
while it is -2.52E-002 if without limiters.

\begin{table}
\centering
\caption{$L^1$ and $L^{\infty}$ errors and orders for Example
\ref{Eg:7} with $u_0(x,y)=\sin^{4}(x+y)$}. \label{tab:6}
\begin{tabular}{lllllllll}
\hline\noalign{\smallskip}
N $\times$ N  &  $L^1$ error & Order&$L^{\infty}$ error& Order & $(\bar u_h)_{min}$ & $(\bar u_h)_{max}$   \\
\noalign{\smallskip} \hline \noalign{\smallskip}
   \multicolumn{7}{c}{\textrm{without limiters}} \\
\noalign{\smallskip} \hline \noalign{\smallskip}
     20  $\times$   20 &    1.70E-03 & &    1.22E-02 & &-2.7162407531E-003 &0.9729124291\\
    40  $\times$   40 &    1.89E-04 &     3.17 &    2.43E-03 &     2.33 &-6.7740627680E-005 &0.9813145773\\
    80  $\times$   80 &    1.97E-05 &     3.26 &    5.51E-04 &     2.14 &-9.6911213559E-006&0.9979503187\\
   160 $\times$   160 &    8.80E-07 &     4.48 &    5.57E-05 &     3.31 &-7.5156683063E-007&0.9993992055\\
   320 $\times$   320 &    4.24E-08 &     4.37 &    3.07E-06 &     4.18 &-8.4287588702E-008&0.9997704325\\
   640 $\times$   640 &    1.35E-09 &     4.97 &    1.25E-07 &     4.62 &-1.7786892420E-008&0.9999547580\\
\noalign{\smallskip} \hline \noalign{\smallskip}
   \multicolumn{7}{c}{\textrm{with limiters}} \\
\noalign{\smallskip} \hline \noalign{\smallskip}
    20  $\times$  20 &    2.37E-03 & &    1.14E-02 & &1.4611105650E-003&  0.9676804528\\
    40  $\times$ 40 &    2.96E-04 &     3.00 &    1.96E-03 &     2.54 &2.3919879074E-004  &0.9819309191\\
    80  $\times$  80 &    2.96E-05 &     3.32 &    5.60E-04 &     1.81 &1.4266765832E-005  &0.9978612888\\
   160 $\times$  160 &    1.27E-06 &     4.56 &    5.57E-05 &     3.33 &1.7175042462E-007 & 0.9993913736\\
   320 $\times$ 320 &    6.26E-08 &     4.33 &    3.07E-06 &     4.18 &1.5078276150E-008 & 0.9997697244\\
   640 $\times$  640 &    2.44E-09 &     4.68 &    1.26E-07 &     4.61 &2.4266981729E-009&  0.9999546678\\
\noalign{\smallskip}\hline
\end{tabular}
\end{table}

\begin{exa}\label{Eg:8}
The two-dimensional inviscid Buckley-Leverett equation with gravitational effects in
the $y$ direction \cite{kurganov2000new} can be written as
\begin{equation}
u_t+f(u)_x+g(u)_y=0,\quad (x,y)\in[-1.5,1.5]\times[-1.5,1.5],
\end{equation}
where
\begin{equation}
f(u)=\frac{u^2}{u^2+(1-u)^2}, \quad g(u)=f(u)(1-5(1-u)^2).
\end{equation}
and the initial condition is
\begin{equation}\label{ic:8}
u(x,y,0)=\left\{\begin{array}{ll}
1, \quad & \textrm{$x^2+y^2<0.5$},\\
0, \quad & \textrm{otherwise}.
\end{array}\right.
\end{equation}

We compute the solution up to time $T=0.5$ with periodic boundary
conditions. The maximum and minimum numerical solutions are listed
in Table \ref{tab:7} for the 5th order FVCW scheme with and without
limiters. Without limiters, we can observe obvious undershoots. They
are completely eliminated by the scheme with limiters. The numerical
solution of FVCW scheme at $T=0.5$ with $N\times N=128\times 128$
are shown in Fig. \ref{Fig:6}, which is similar to the result in
\cite{kurganov2000new}. The contour plot in Fig. \ref{Fig:6} shows
that the present scheme produces high resolution numerical solutions
without significant spurious oscillations.

\begin{table}
\centering
\caption{Maximum and minimum numerical solutions for Example
\ref{Eg:8} at $T = 0.5$.} \label{tab:7}
\begin{tabular}{c|cc|cc}
\hline%\noalign{\smallskip}
&  with limiters    &               & without limiters \\
\hline%\noalign{\smallskip}
N $\times$ N  &  $(\bar u_h)_{min}$ & $(\bar u_h)_{max}$ &$(\bar u_h)_{min}$ & $(\bar u_h)_{max}$ \\
%\noalign{\smallskip}
\hline %\noalign{\smallskip}
    8 $\times$ 8 &   5.2763420443E-004 & 0.9182487023   &-5.3290980853E-002  &1.0218414115     \\
    16 $\times$ 16 & 4.9239390613E-007 & 0.9514705351   &-2.5761986049E-002  &0.9554234787 \\
    32 $\times$ 32 & 7.8904099515E-012 & 0.9840518251   &-4.7749433037E-003  &0.9841847426\\
    64 $\times$ 64 & 8.6575061287E-017 & 0.9971065091   &-2.8304959899E-007  &0.9971278020 \\
   128 $\times$128 & 1.3048650350E-025 & 0.9998762415   &-1.6326732186E-008  &0.9998768328\\
   128 $\times$128 & 5.3604284574E-041 & 0.9999998498   &-1.8832845160E-008  &0.9999998751\\
\hline
\end{tabular}
\end{table}

\begin{figure}[htp]
\subfigure{\includegraphics[width=0.5\textwidth]{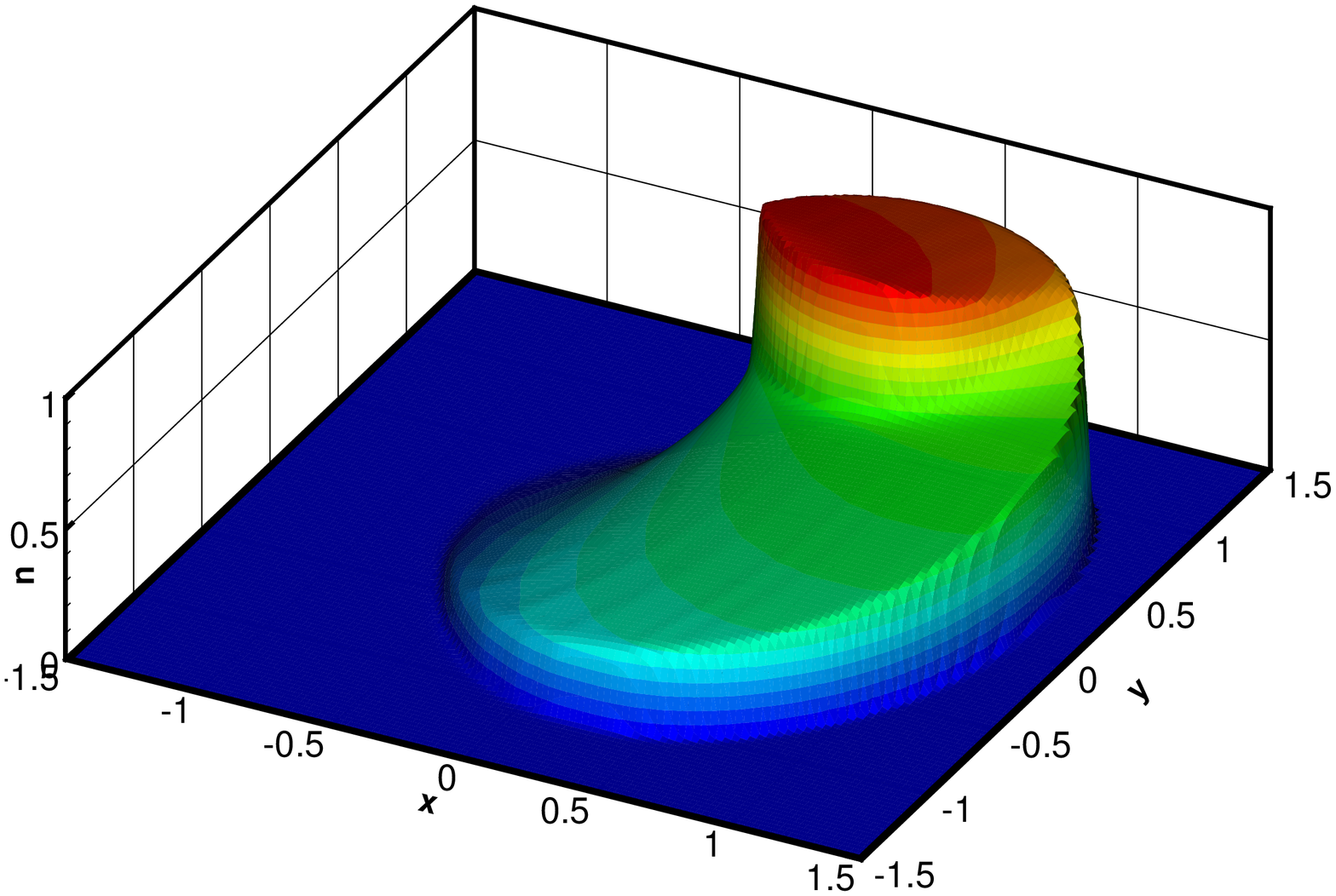}}
\subfigure{\includegraphics[width=0.5\textwidth]{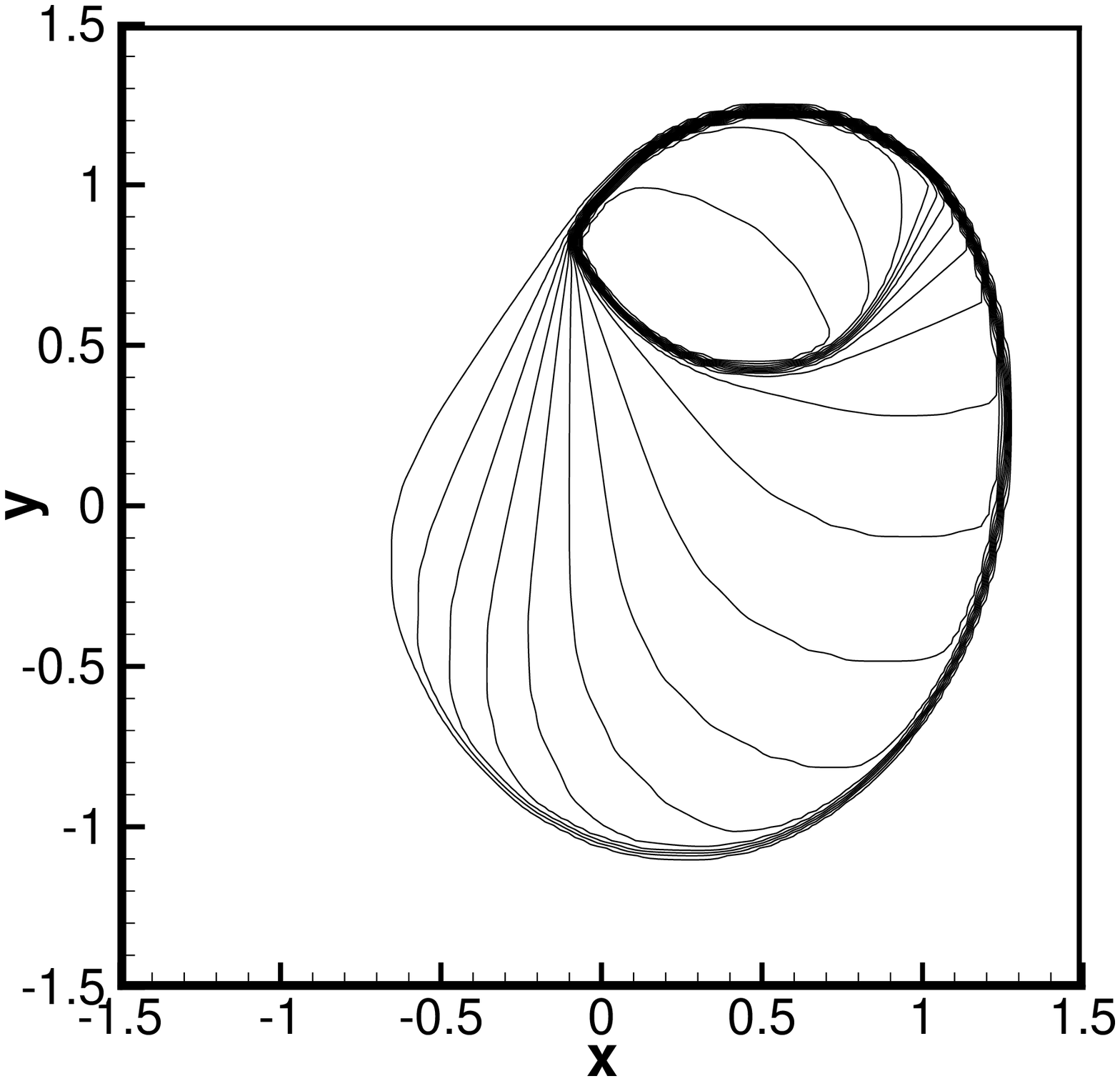}}
\caption{Numerical solution for Example \ref{Eg:8} with FVCW scheme
at T=0.5 on a $128\times 128$ mesh.}\label{Fig:6}
\end{figure}
\end{exa}

\subsection{Incompressible flow}
In this section, both the FVCW scheme and the finite volume compact (FVC)
scheme (which is the FVCW scheme with optimal linear weights)
are considered. The FVCW scheme with limiters is denoted as ``FVCW-W''
and it is denoted as ``FVCW-O'' if without limiters. Similarly,
the FVC scheme with and without limiters are denoted as ``FVC-W'' and ``FVC-O''.

\begin{figure}[htp]
\subfigure{\includegraphics[width=0.5\textwidth]{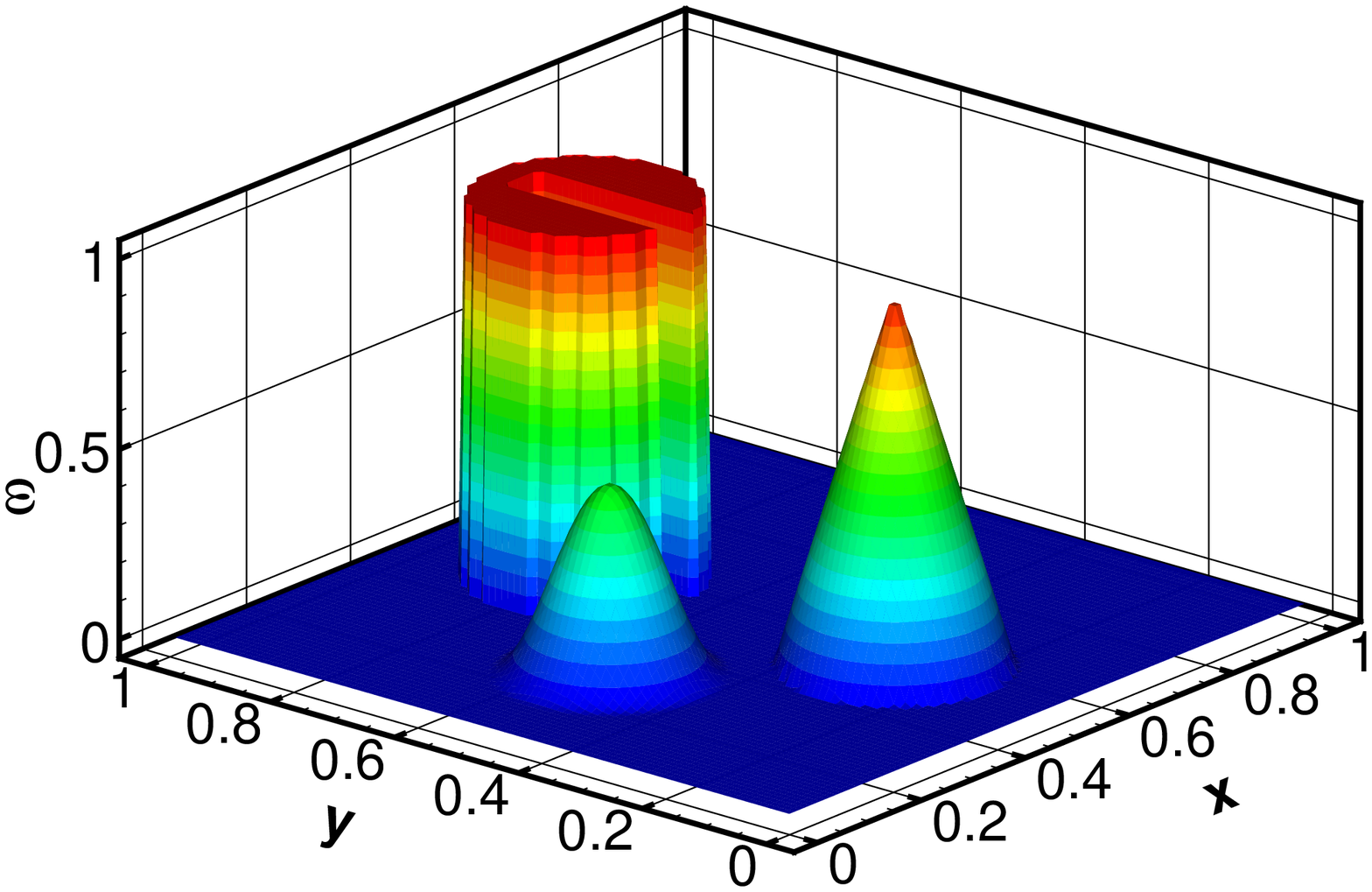}}
\subfigure{\includegraphics[width=0.5\textwidth]{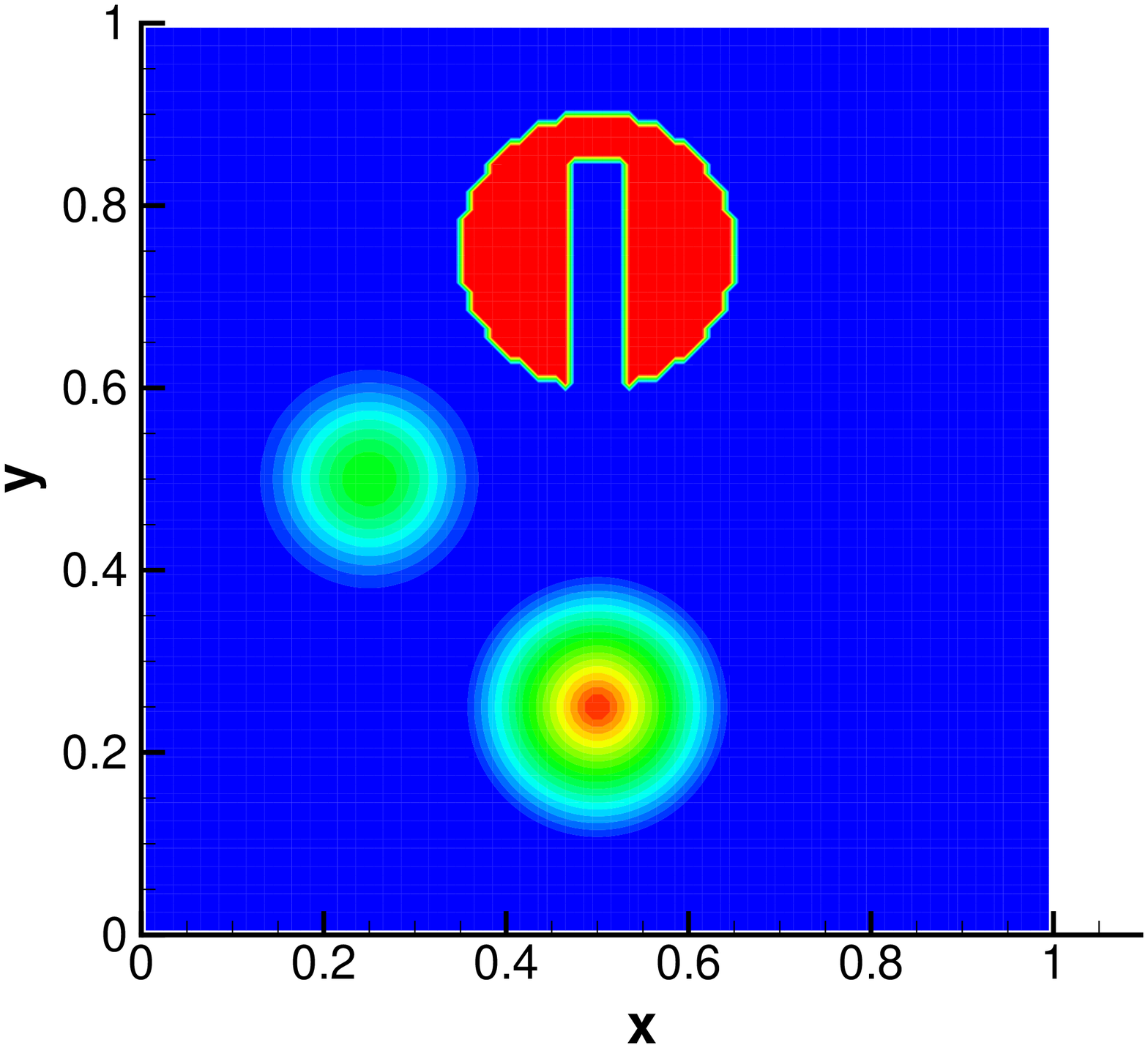}}
\caption{Initial profile.}\label{Fig:7}
\end{figure}

\begin{figure}[htp]
\subfigure{\includegraphics[width=0.5\textwidth]{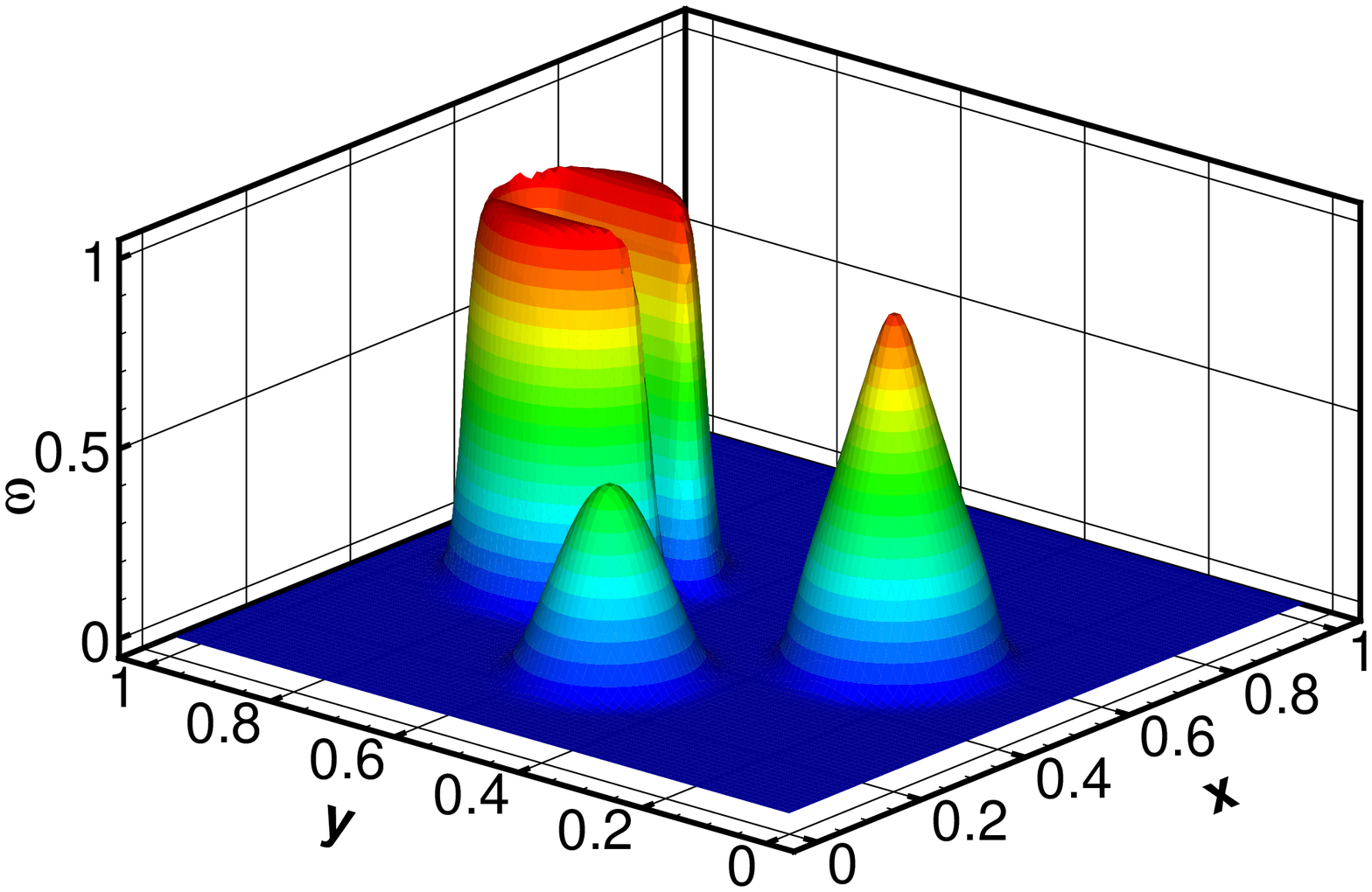}}
\subfigure{\includegraphics[width=0.5\textwidth]{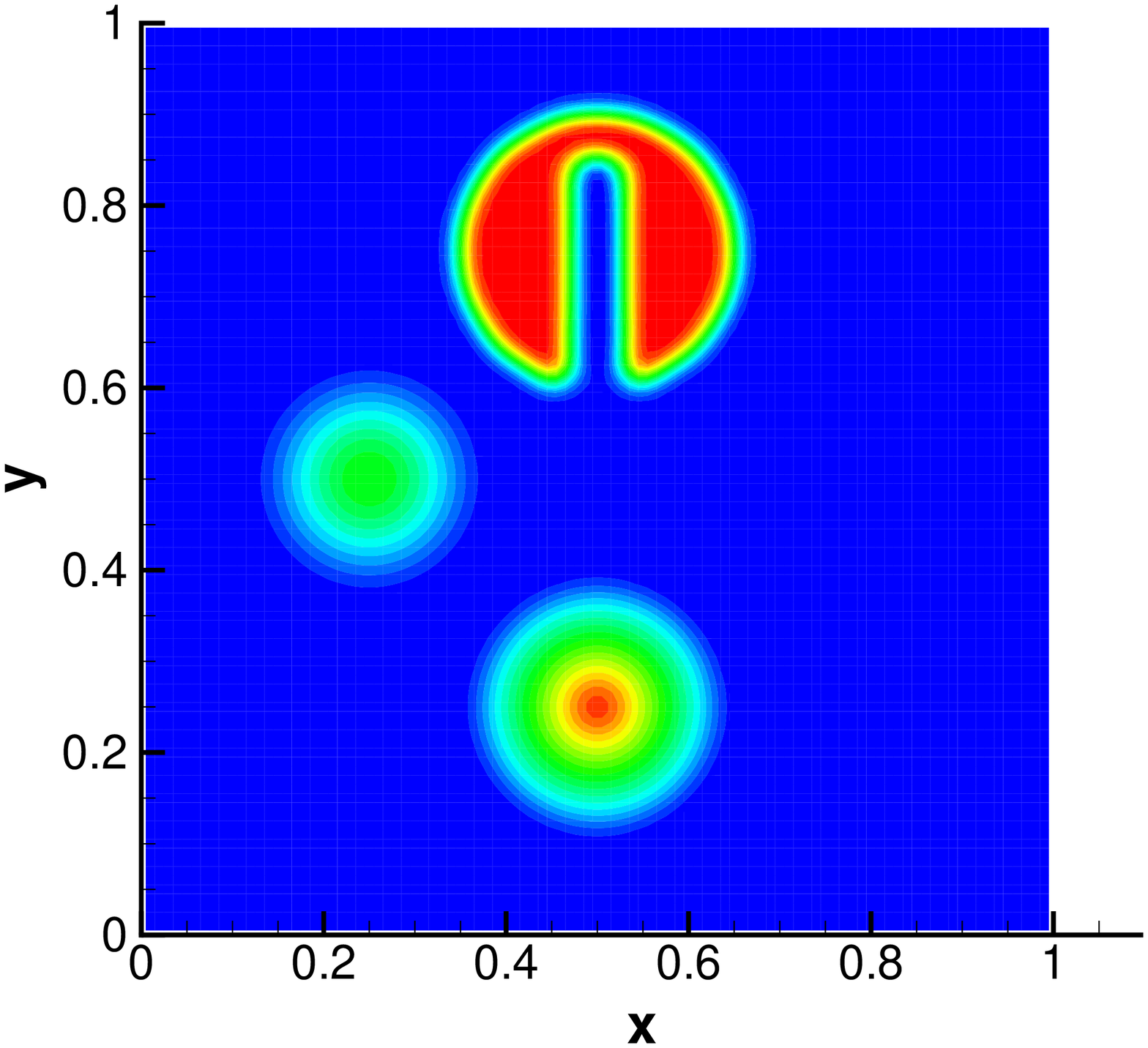}}
\caption{Numerical results for Example \ref{Eg:9} at $T=2\pi$ with
limiters.}\label{Fig:8}
\end{figure}

\begin{figure}[htp]
\subfigure[$y=0.25$]{\includegraphics[width=0.5\textwidth]{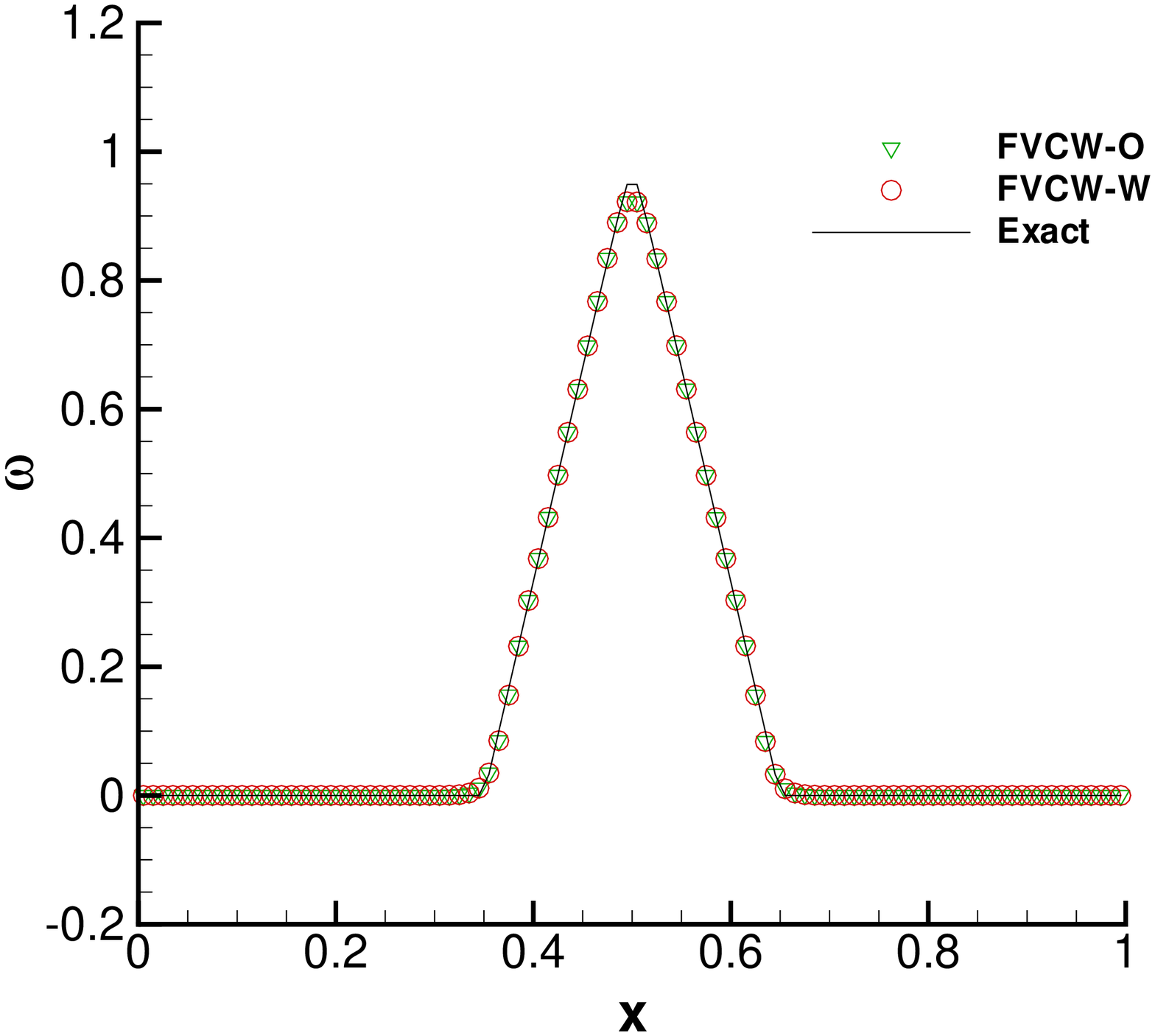}}
\subfigure[$y=0.25$]{\includegraphics[width=0.5\textwidth]{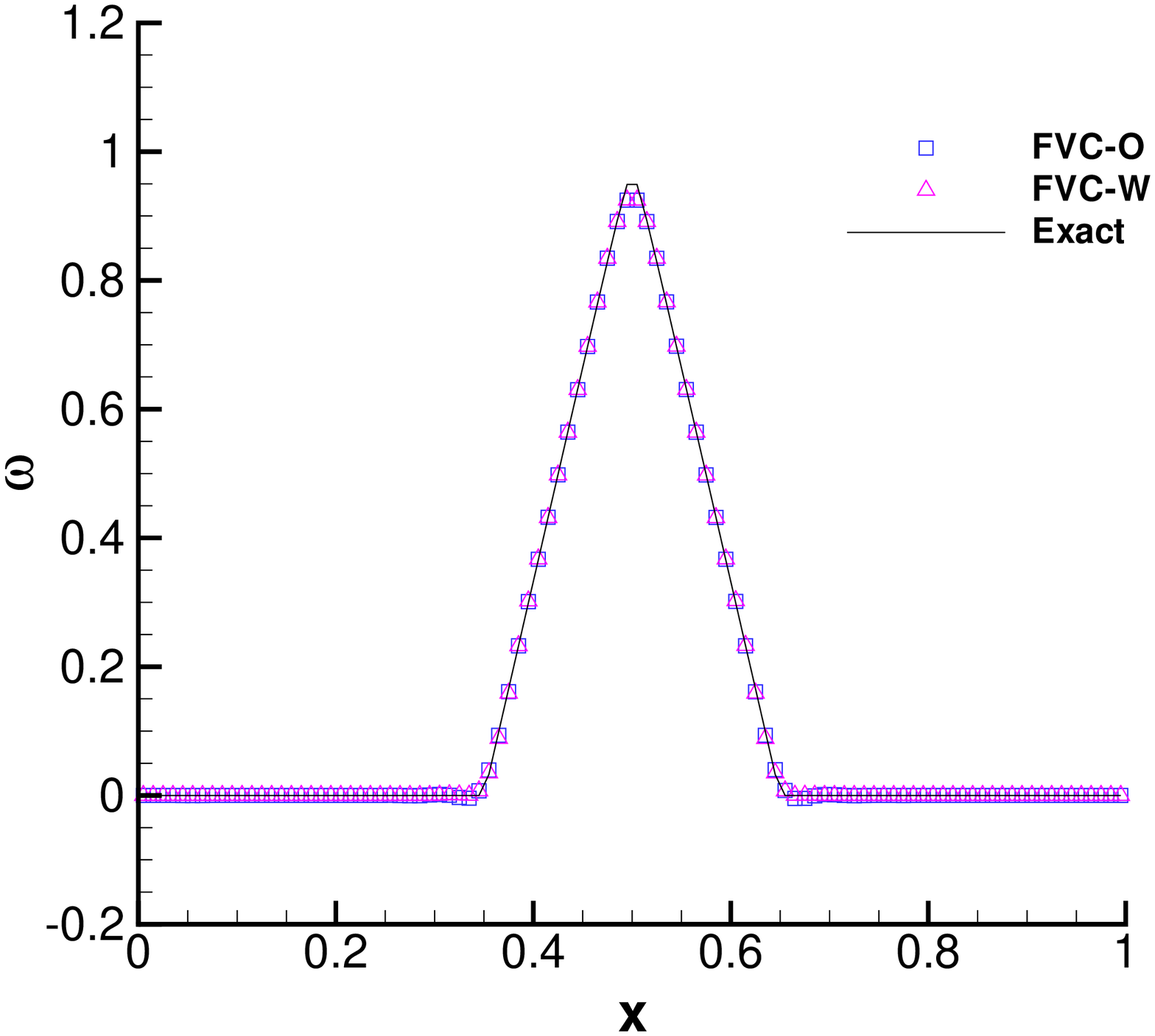}}
\subfigure[$y=0.75$]{\includegraphics[width=0.5\textwidth]{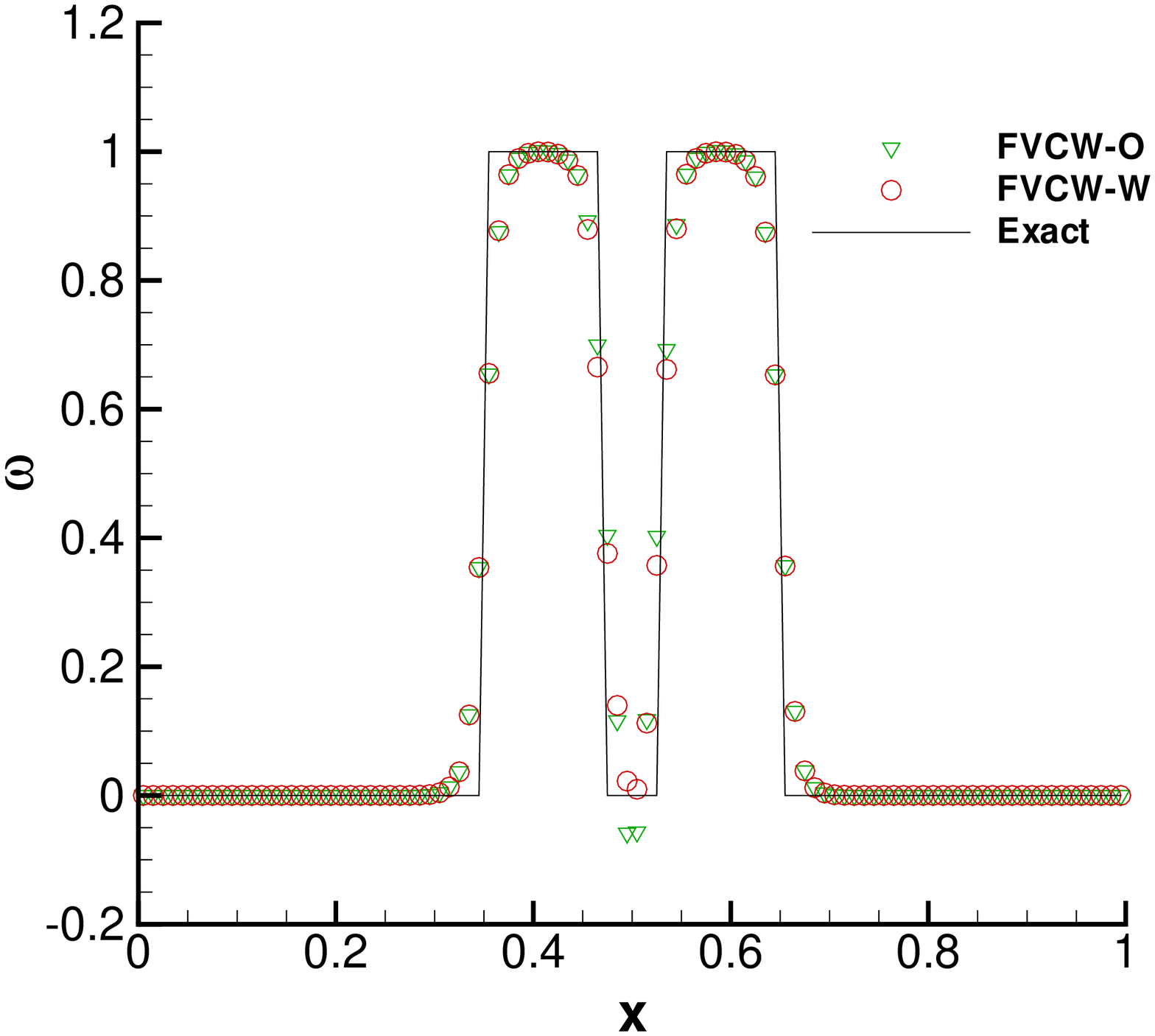}}
\subfigure[$y=0.75$]{\includegraphics[width=0.5\textwidth]{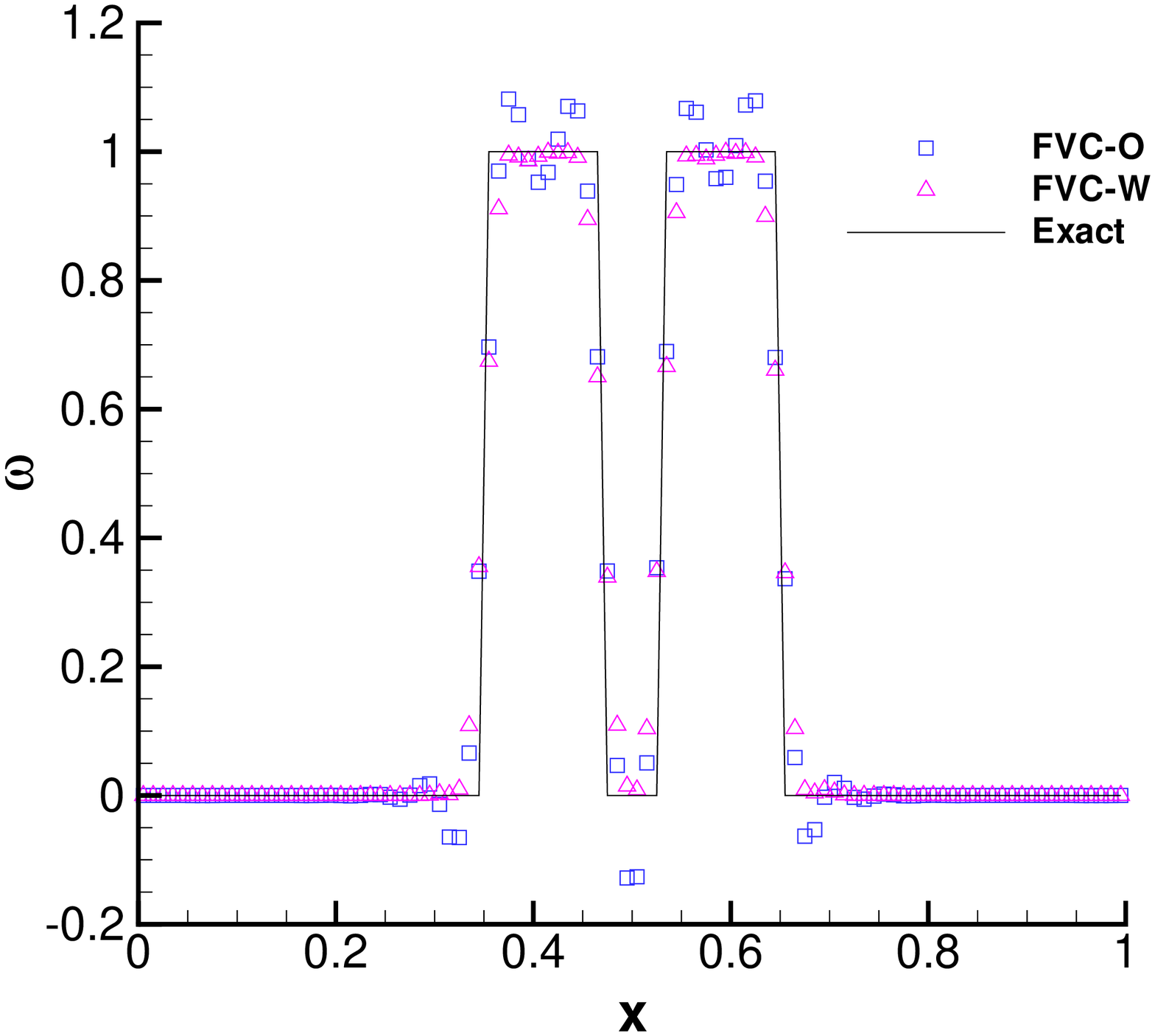}}
\subfigure[$x=0.5$]{\includegraphics[width=0.5\textwidth]{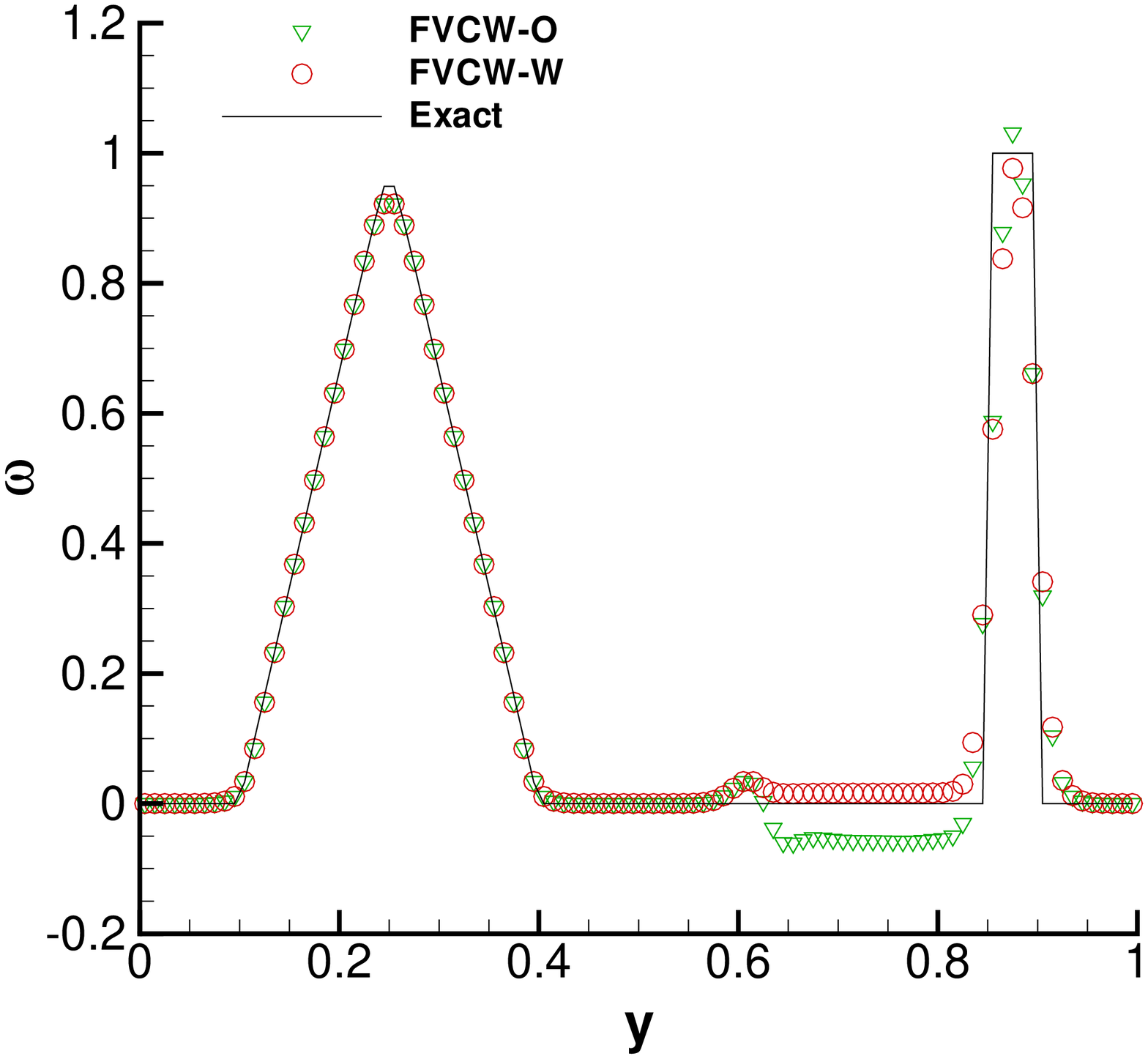}}
\subfigure[$x=0.5$]{\includegraphics[width=0.5\textwidth]{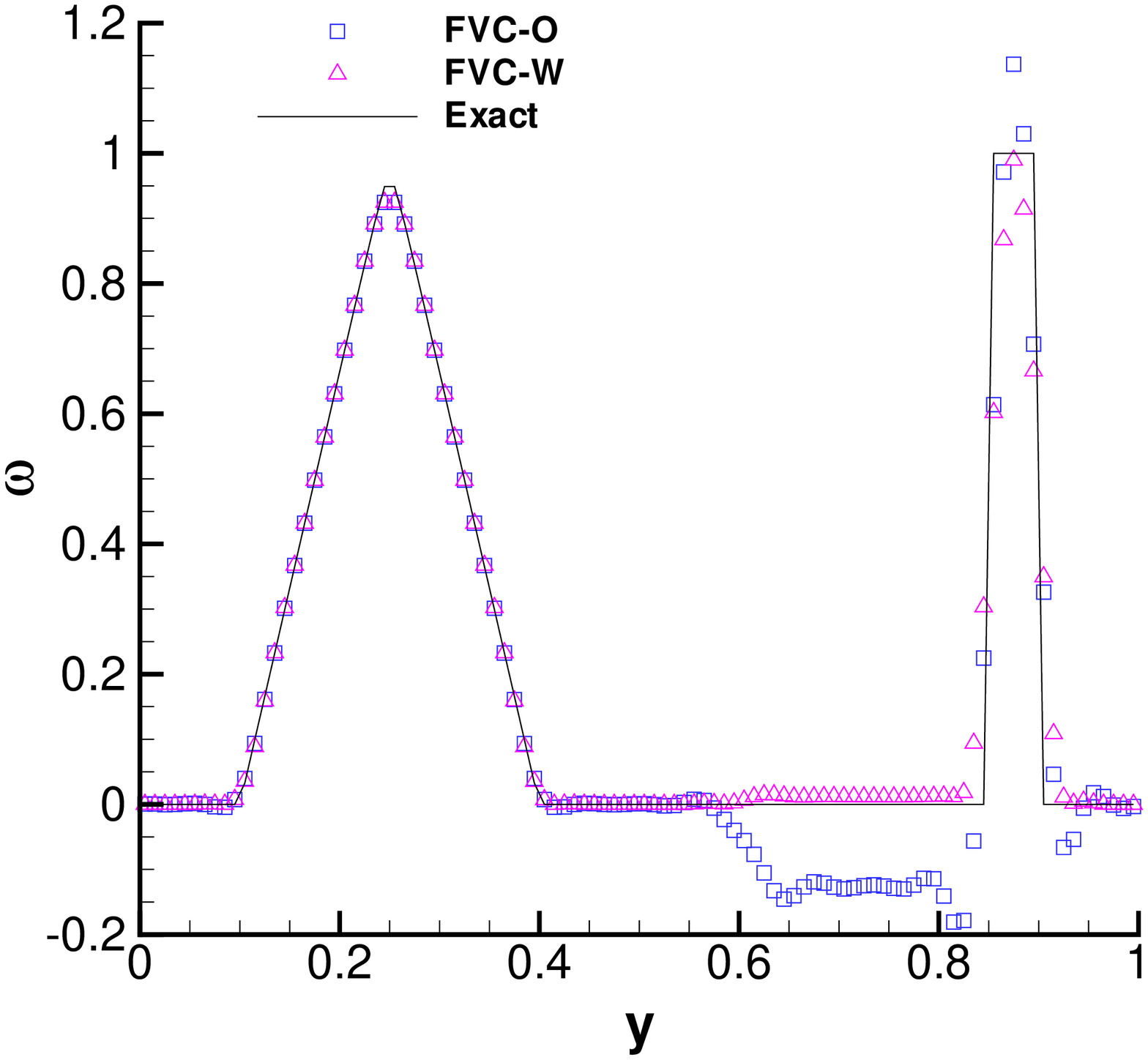}}
\caption{Cross-sections of the numerical solutions with the exact
solutions for Example \ref{Eg:9} at $T=2\pi$. Left: FVCW scheme;
Right: FVC scheme.}\label{Fig:9}
\end{figure}

\begin{exa}\label{Eg:9}
(Rigid body rotation) Consider the following rigid body rotation problem \cite{leveque1996high}
\begin{equation}\label{sbr}
\omega_t+(-(y-0.5)\omega)_x+((x-0.5)\omega)_y=0, \quad (x,y)\in[0,1]\times[0,1].
\end{equation}
\end{exa}
The initial profile consists of a smooth hump, a cone and a slotted
cylinder \cite{leveque1996high,park2010multi} as in Fig.
\ref{Fig:7}. The numerical solution of the FVCW scheme with limiters
after one period of evolution on a mesh of $N\times N =100\times
100$ is shown in Fig. \ref{Fig:8}. The numerical solutions of the
FVCW scheme and the FVC scheme with and without limiters at $T=2\pi$
are plotted in Fig. \ref{Fig:9}, by comparing with the exact
solutions. With limiters, the minimum numerical value of the FVCW
scheme is 4.46E-015 and the maximum value is 0.99984. However, they
are -6.24E-002 and 1.05022 if without limiters. The numerical
results are similar to those in \cite{qiu2011conservative}.

\begin{exa}\label{Eg:10}
(Swirling deformation flow) The swirling deformation flow \cite{leveque1996high} is
\begin{equation}
\omega_t+(\sin^2{(\pi x)}\sin{(2\pi y)}g(t)\omega)_x+(-\sin^2{(\pi
y)}\sin{(2\pi x)}g(t)\omega)_y=0, \nonumber \\
(x,y)\in[0,1]\times[0,1],
\label{sdf}
\end{equation}
with $g(t)=\cos{(\pi t/T)}$ on the time interval $0\leq t \leq T$.
The initial condition is the same as in Fig. \ref{Fig:7}.

The flow slows down and reverses direction when the initial data is
recovered at time $T$. We take $T=1.5$. The numerical solution of
the FVCW scheme on a mesh of $N \times N=100\times 100$ at $T=1.5$
is shown in Fig. \ref{Fig:10}. Similarly the cuts of the numerical
solution at $T$ are plotted in Fig. \ref{Fig:11} by comparing with
the exact solutions. With limiters, the minimum numerical value is
3.50E-018 and the maximum value is 0.99967. However, they are
-6.76E-002 and 1.09007 if without limiters. At time $T/2$, the
solution is quite deformed from the initial data, see Fig.
\ref{Fig:12}. The minimum and maximum values with limiters are
5.91E-019 and 0.99903. Without limiters, they are -6.60E-002 and
1.06939 respectively.

\begin{figure}[htp]
\subfigure{\includegraphics[width=0.5\textwidth]{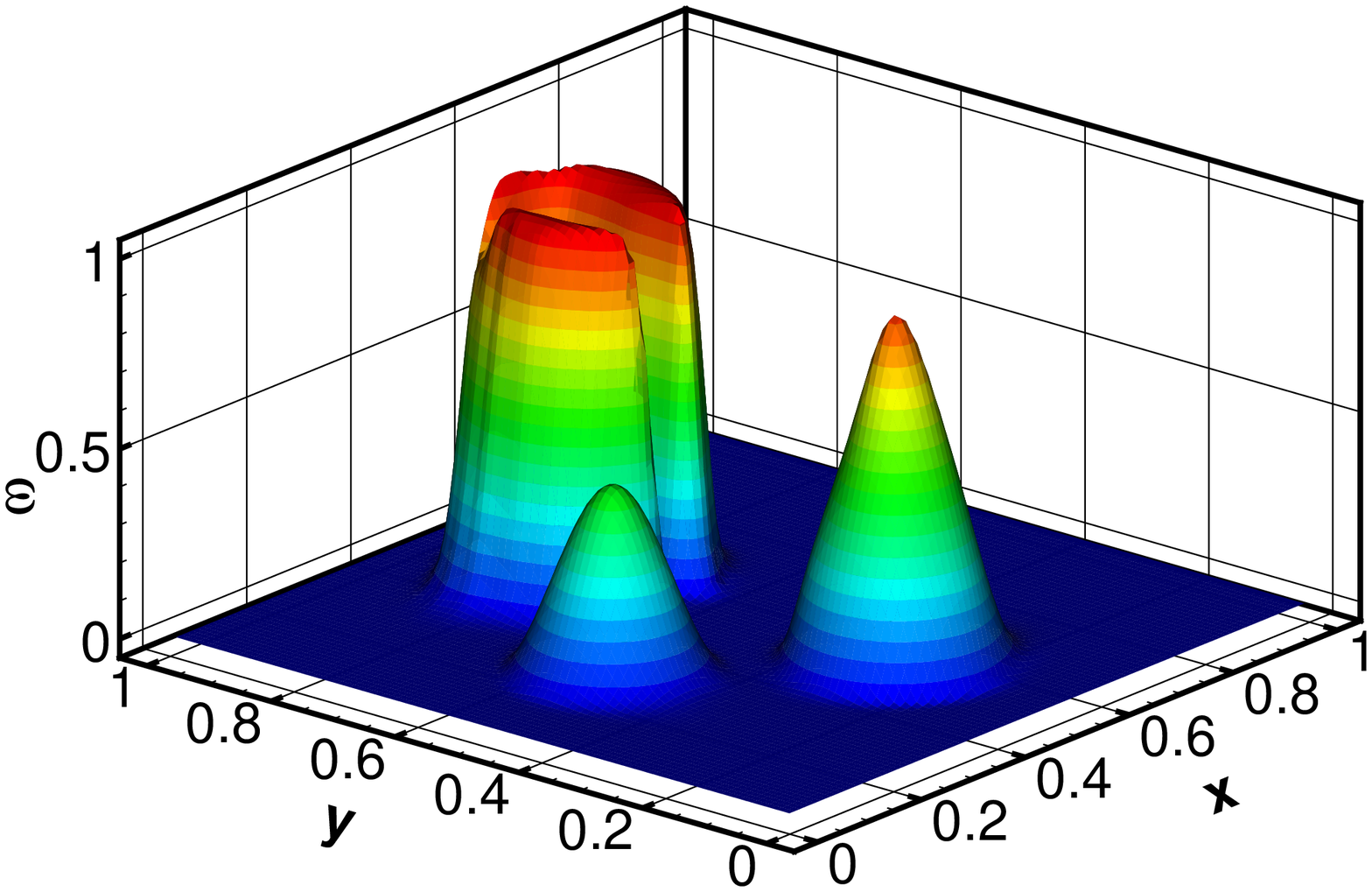}}
\subfigure{\includegraphics[width=0.5\textwidth]{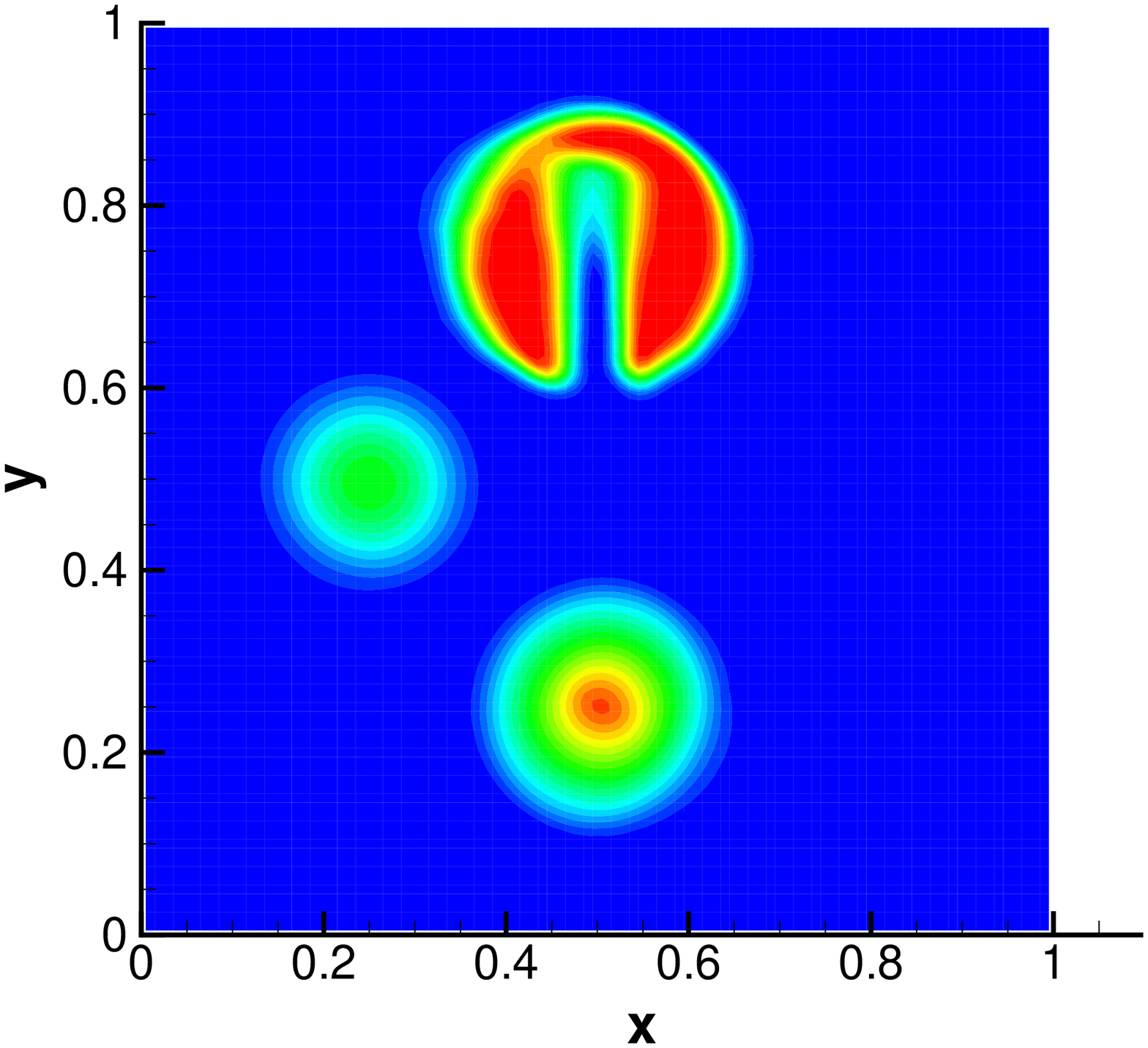}}
\caption{Numerical results for Example \ref{Eg:10} at
$T=1.5$.}\label{Fig:10}.
\end{figure}

\begin{figure}[htp]
\subfigure[$y=0.25$]{\includegraphics[width=0.5\textwidth]{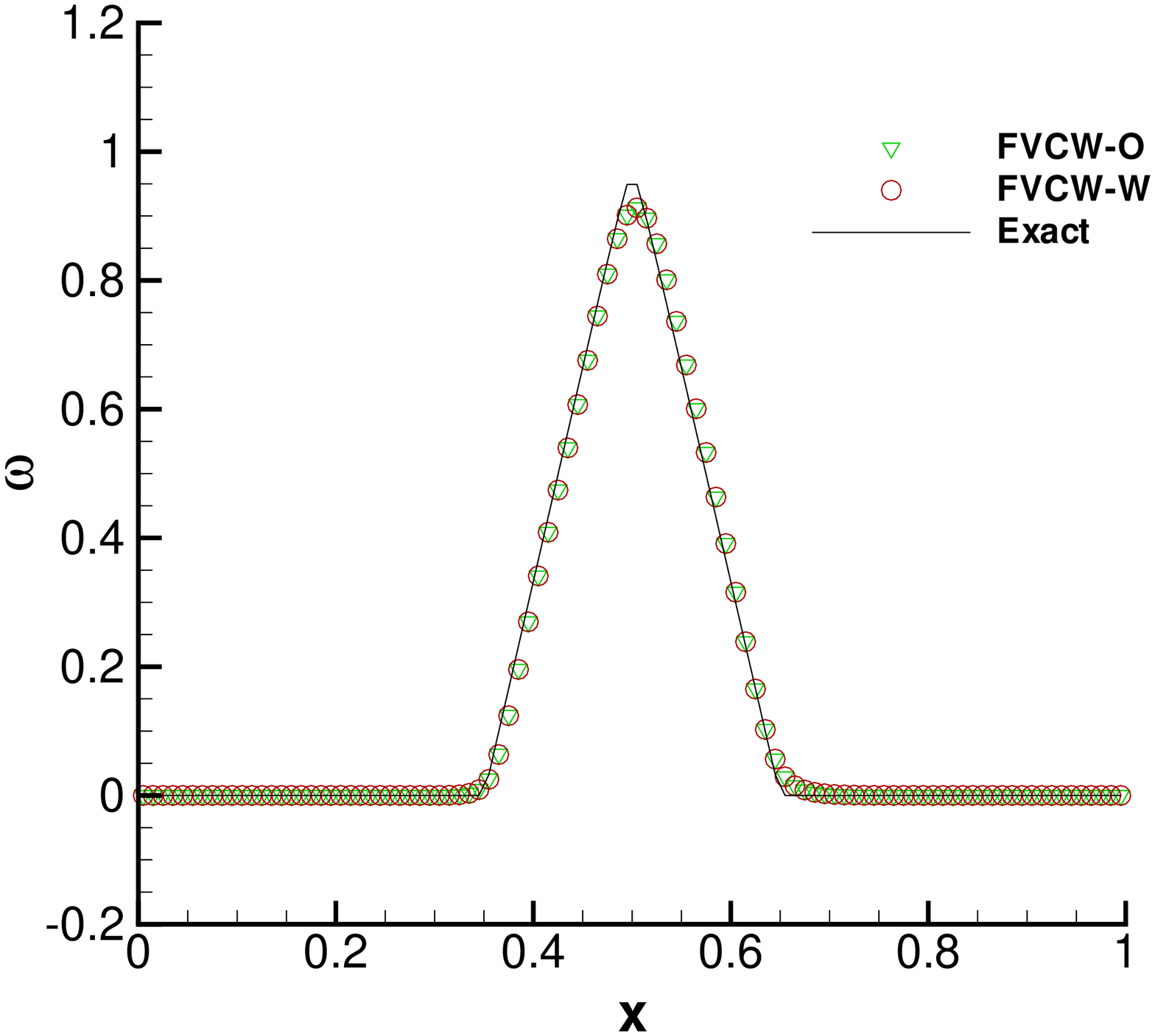}}
\subfigure[$y=0.25$]{\includegraphics[width=0.5\textwidth]{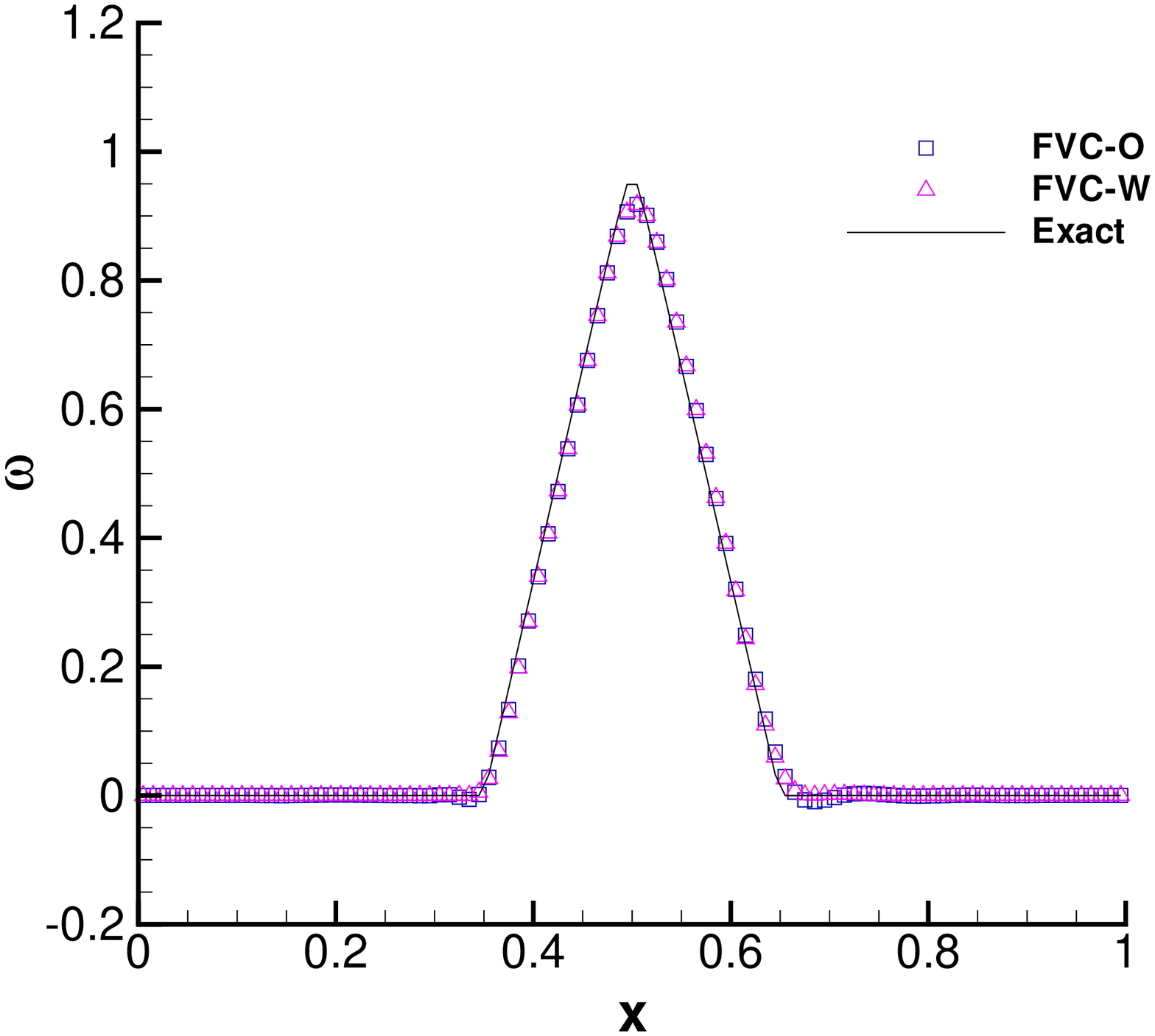}}
\subfigure[$y=0.75$]{\includegraphics[width=0.5\textwidth]{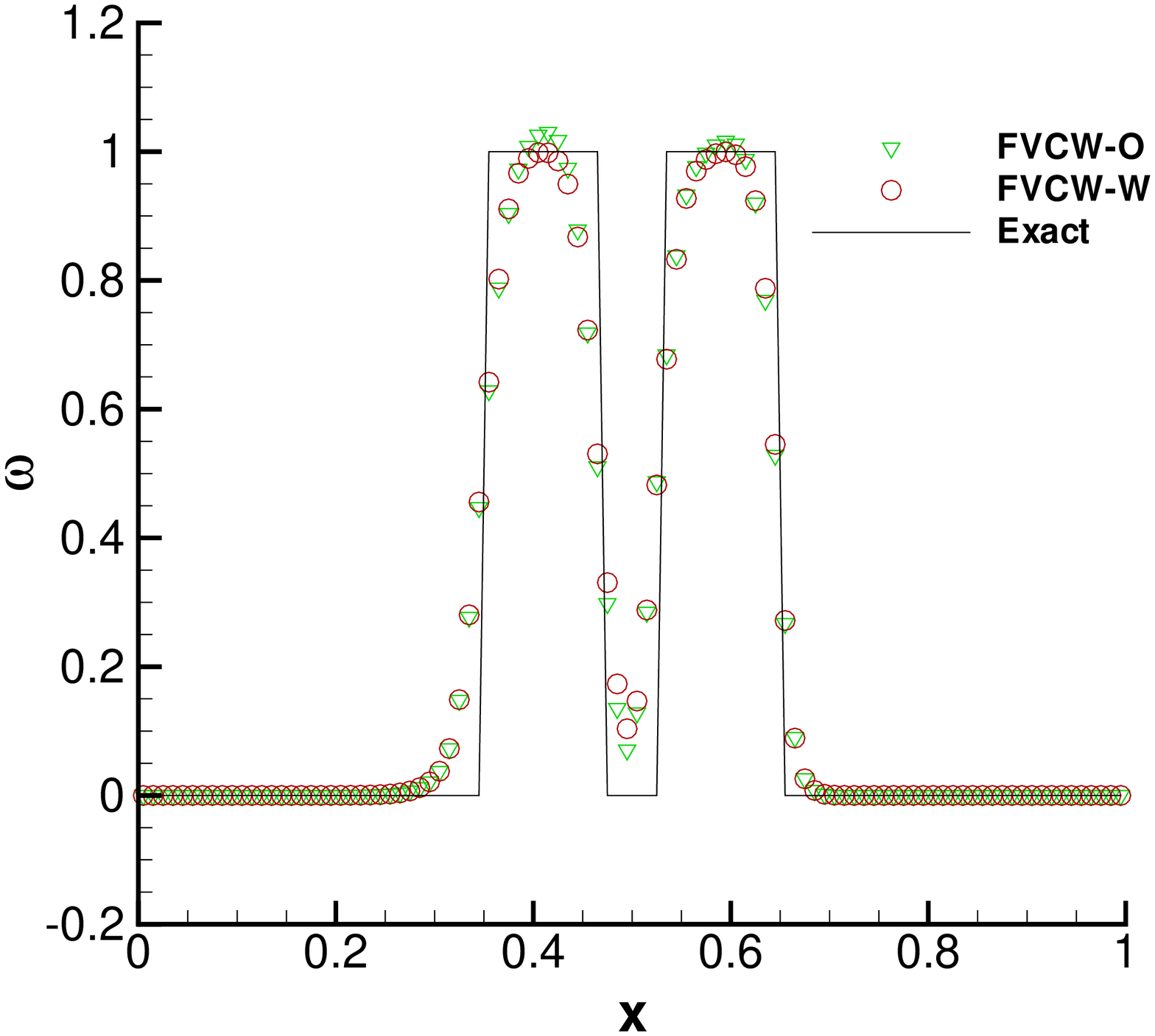}}
\subfigure[$y=0.75$]{\includegraphics[width=0.5\textwidth]{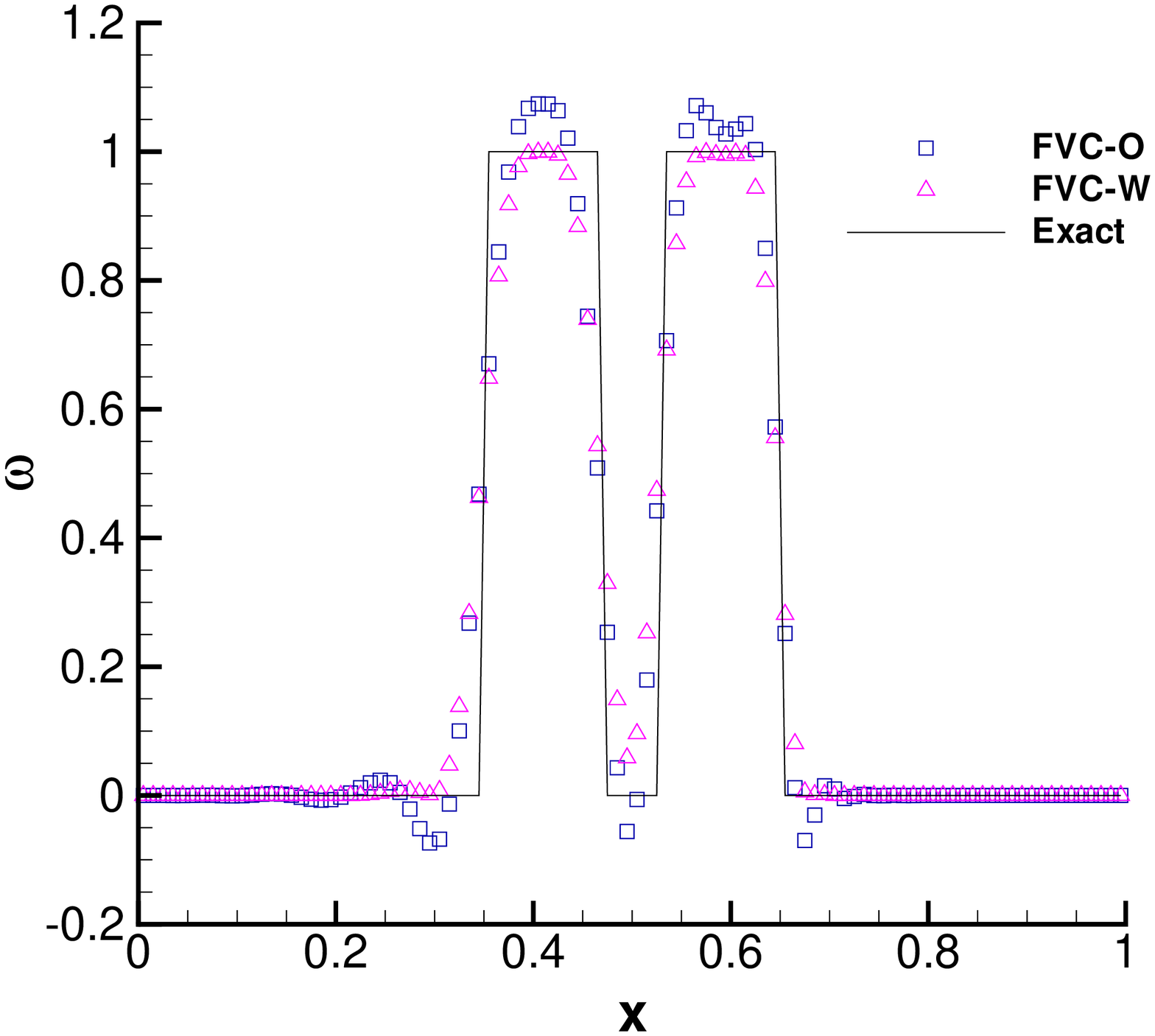}}
\subfigure[$x=0.5$]{\includegraphics[width=0.5\textwidth]{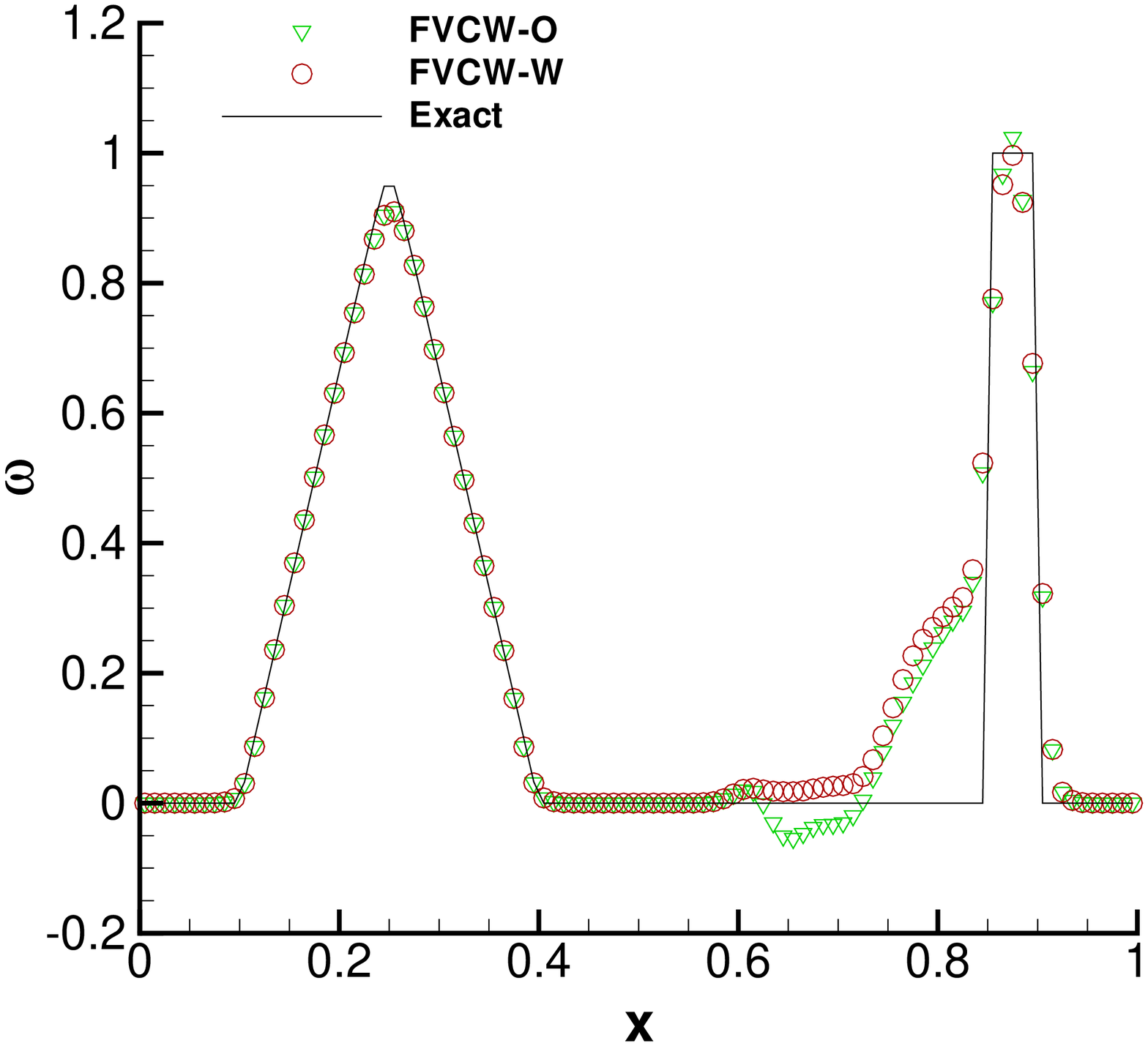}}
\subfigure[$x=0.5$]{\includegraphics[width=0.5\textwidth]{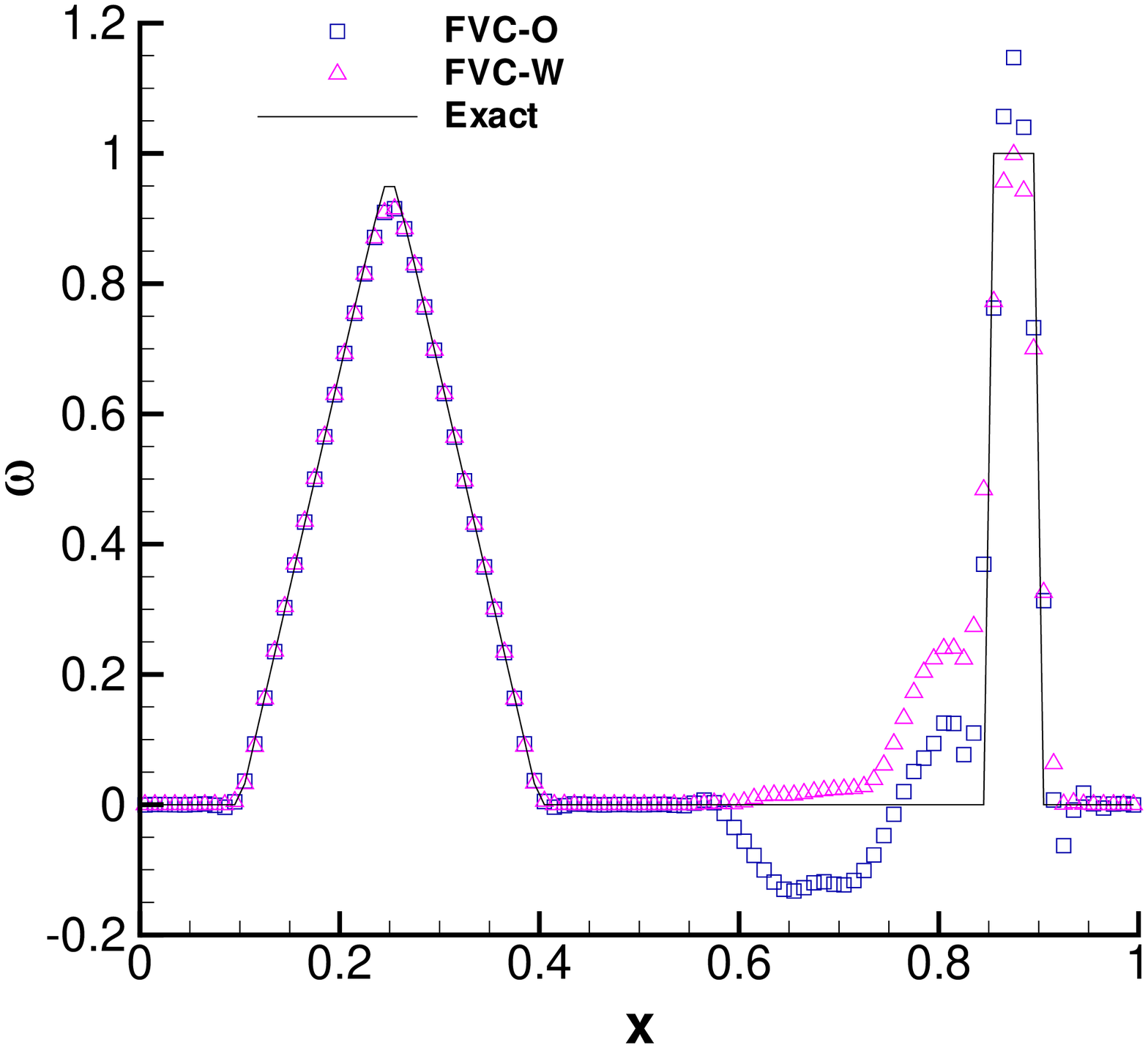}}
\caption{Cross-sections of the numerical solutions with the exact
solutions for Example \ref{Eg:10} at $T=1.5$.Left: FVCW scheme;
Right: FVC scheme.}\label{Fig:11}
\end{figure}

\begin{figure}[htp]
\subfigure{\includegraphics[width=0.5\textwidth]{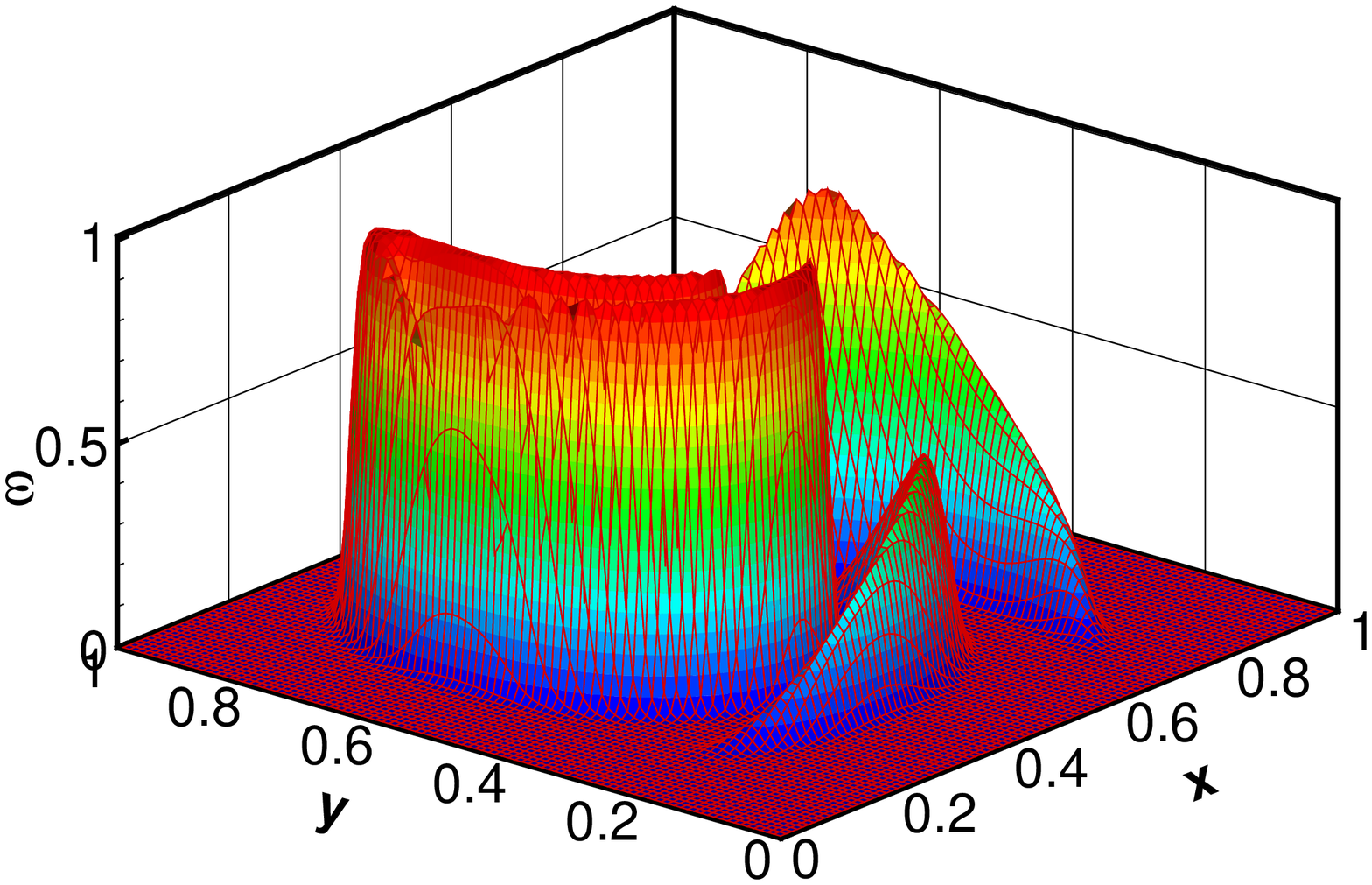}}
\subfigure{\includegraphics[width=0.5\textwidth]{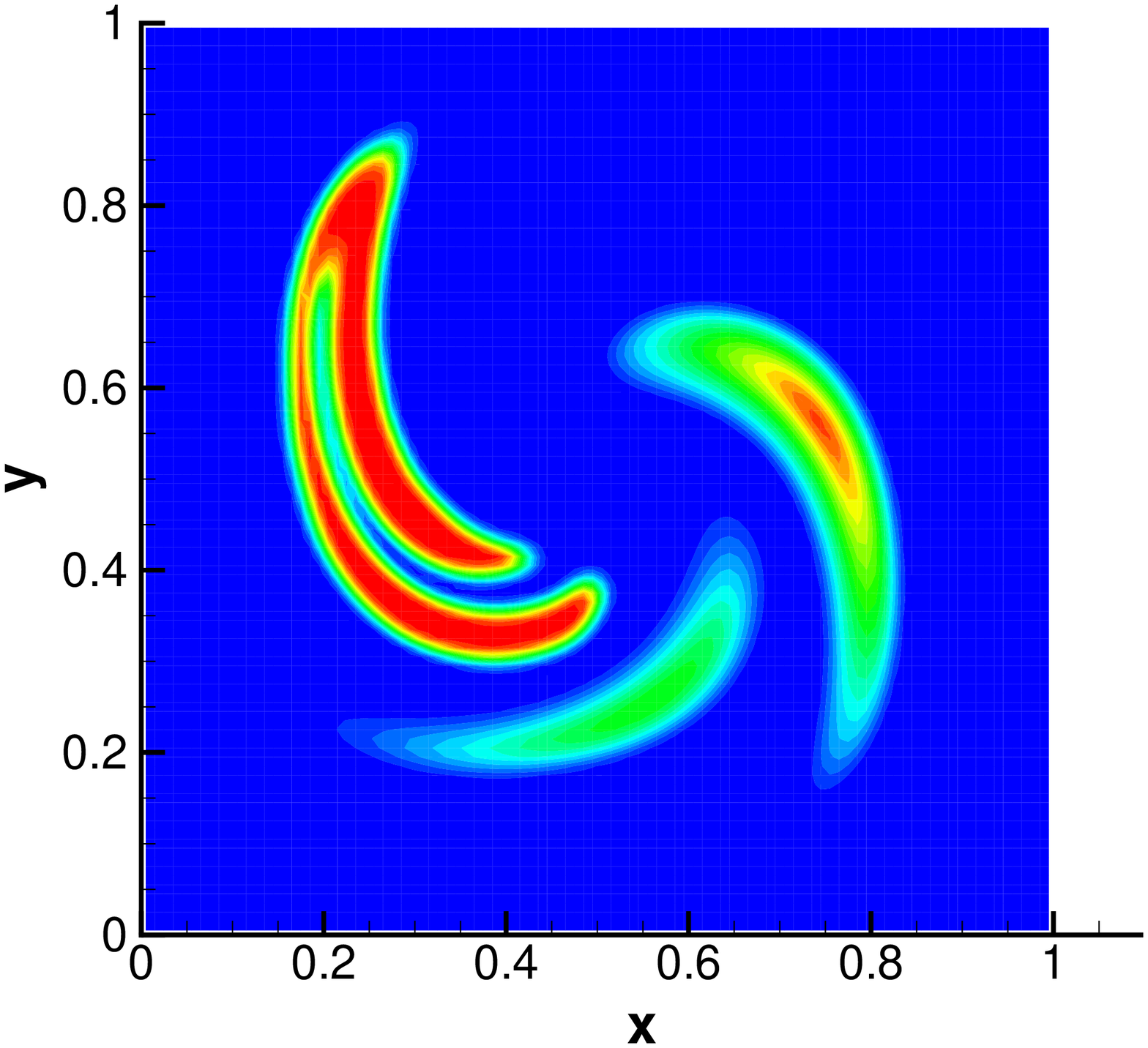}}
\caption{The deformed numerical results for Example \ref{Eg:10} at
$T/2$.}\label{Fig:12}
\end{figure}
\end{exa}

\begin{exa}\label{Eg:11}
 In this example, we consider two-dimensional incompressible equations \eqref{incom}.
 The computational domain is $[0, 2\pi]\times [0, 2\pi]$ with
  periodic boundary conditions.
For this problem with the smooth exact solution
$\omega(x,y,t)=-2\sin(x)\sin(y)$, at $T=1$ the designed 5th order of accuracy
can be clearly observed in Table \ref{tab:8}, and the numerical
solutions with limiters are all with in the range [-2,2].
\begin{table}
\centering
\caption{$L^1$ and $L^{\infty}$ errors and orders for Example
\ref{Eg:11} with $\omega_0(x,y,0)=-2\sin(x)\sin(y)$, with limiters.}
\label{tab:8}
\begin{tabular}{lllllllll}
\hline\noalign{\smallskip}
N $\times$ N  &  $L^1$ error & Order &$L^{\infty}$ error& Order & $(\bar \omega_h)_{min}$ & $(\bar \omega_h)_{max}$   \\
\noalign{\smallskip} \hline \noalign{\smallskip}
    20  $\times$  20 &     1.920E-06 & &    8.349E-06 & &-1.9350632621        &  1.9350632621 \\
    40  $\times$  40 &     3.769E-08 &     5.67 &    2.706E-07 &     4.95 & -1.9836046769     &   1.9836046769\\
    80  $\times$  80 &     7.533E-10 &     5.64 &    8.674E-09 &     4.96 & -1.9958910455   &     1.9958910455     \\
   160  $\times$ 160 &     1.669E-11 &     5.50 &    2.724E-10 &     4.99 & -1.9989721276    &    1.9989721276 \\
    320  $\times$ 320 &     4.170E-13 &     5.32 &    8.524E-12 &     5.00 &
-1.9997429923      &  1.9997429923     \\
\noalign{\smallskip}\hline
\end{tabular}
\end{table}
\end{exa}

\begin{exa}\label{Eg:12}
(The double shear layer problem) We then consider (\ref{incom}) for the double shear layer problem on the domain $[0, 2\pi]\times [0, 2\pi]$ with the initial conditions given by
\begin{equation}\label{ic:12}
\omega(x,y,0)=\left\{\begin{array}{ll}
\delta \cos(x)-\frac{1}{\rho}sech^2(\frac{1}{\rho}(y-\frac{\pi}{2})), \quad & \textrm{$y\leq \pi$},\\
\delta
\cos(x)+\frac{1}{\rho}sech^2(\frac{1}{\rho}(\frac{3\pi}{2}-y)), \quad &
\textrm{otherwise},
\end{array}\right.
\end{equation}
where $\rho=\frac{\pi}{15}$ and $\delta=0.05$.

The contours of the vorticity at $T=6$ with meshes of $64\times 64$
and $128 \times 128$ for the FVCW scheme with and without limiters
are shown in Fig. \ref{Fig:13}. The contours of the vorticity at
$T=8$ are shown in Fig. \ref{Fig:14}. For this problem, we can not
see any visible difference of the results between with and without
limiters, but the numerical solution with limiters is within the
range $[-\delta-\frac{1}{\rho}, \delta+\frac{1}{\rho}]$. The results
are similar to those in
\cite{zhang2010maximum,liu2000high}.

\begin{figure}[htp]
\subfigure[$128\times 128$ mesh, FVCW with
limiter]{\includegraphics[width=0.5\textwidth]{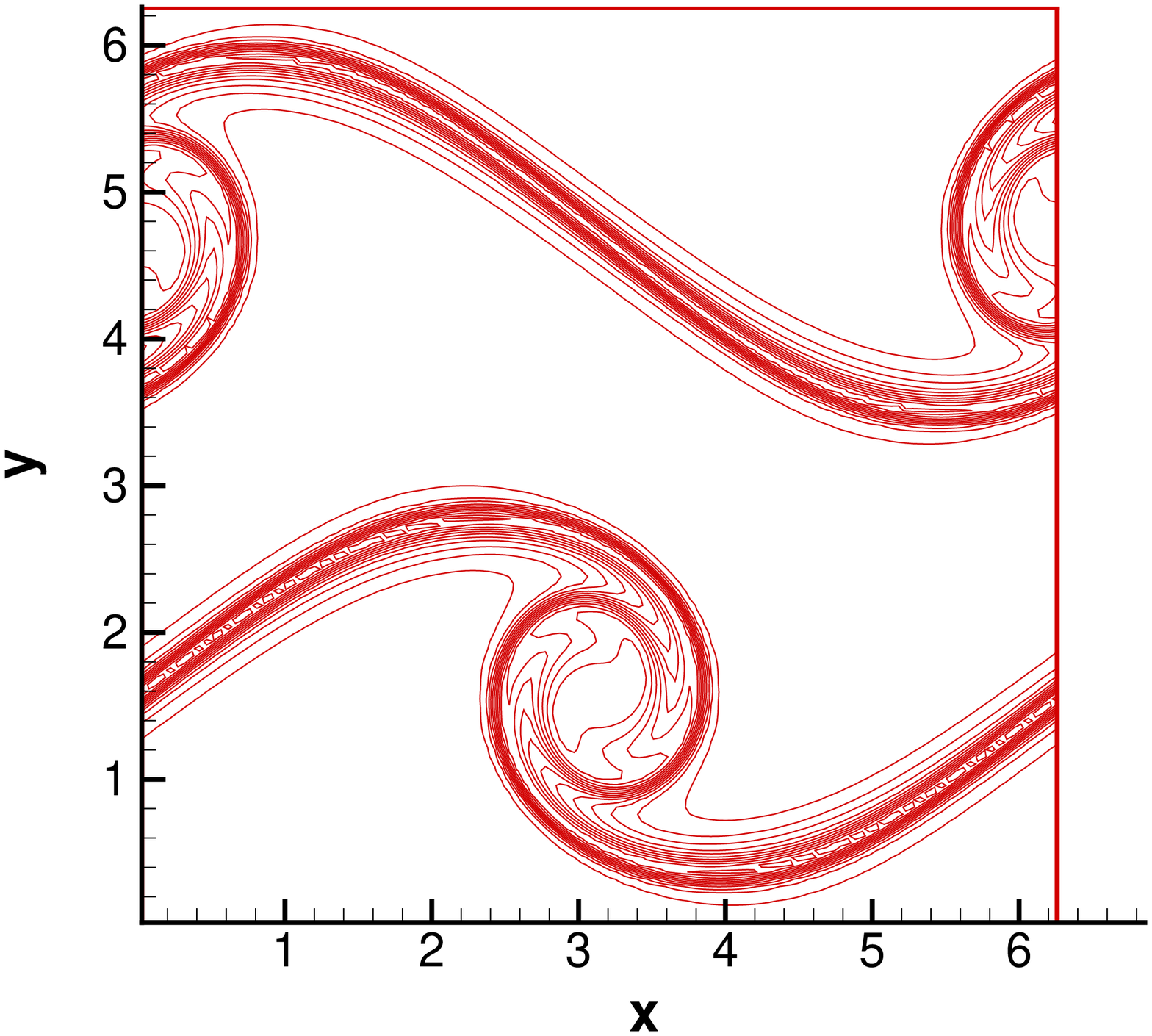}}
\subfigure[$128\times 128$ mesh, FVCW without
limiter]{\includegraphics[width=0.5\textwidth]{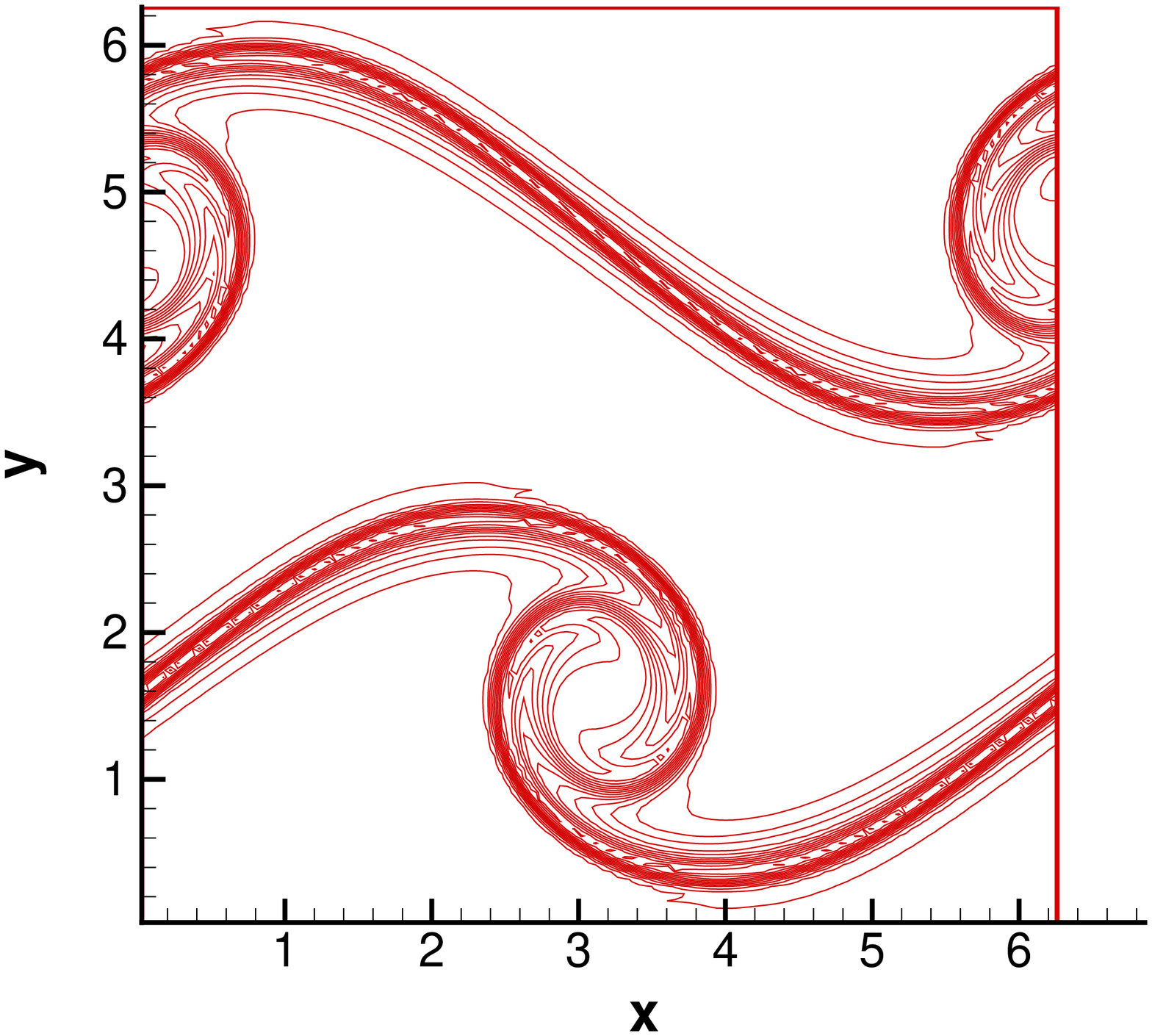}}
\caption{Numerical results for Example \ref{Eg:12} at $T=6$,
$64\times 64$(left),$128\times 128$(right).}\label{Fig:13}.
\end{figure}

\begin{figure}[htp]
\subfigure[$128\times 128$ mesh, FVCW with
limiter]{\includegraphics[width=0.5\textwidth]{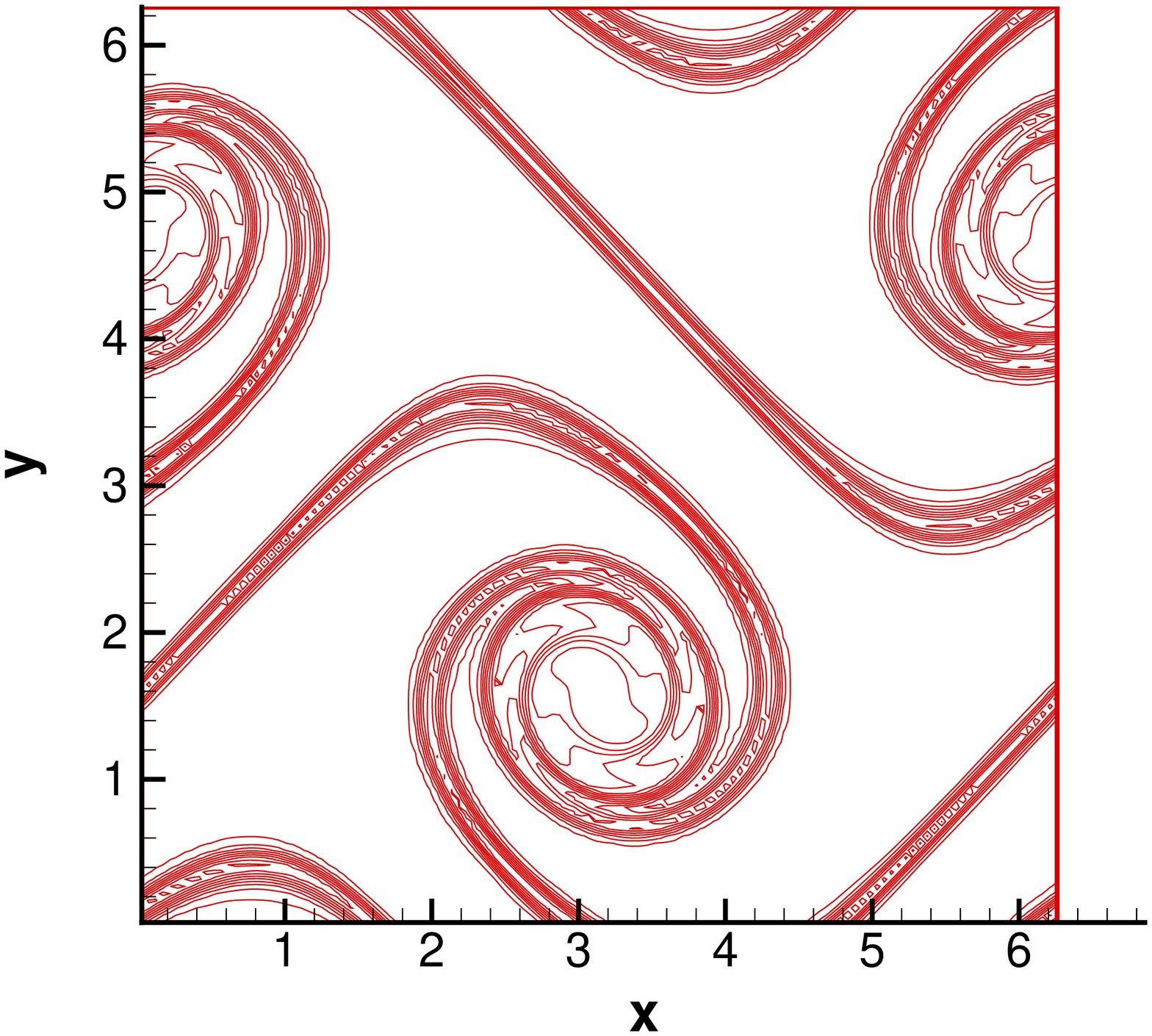}}
\subfigure[$128\times 128$ mesh, FVCW without
limiter]{\includegraphics[width=0.5\textwidth]{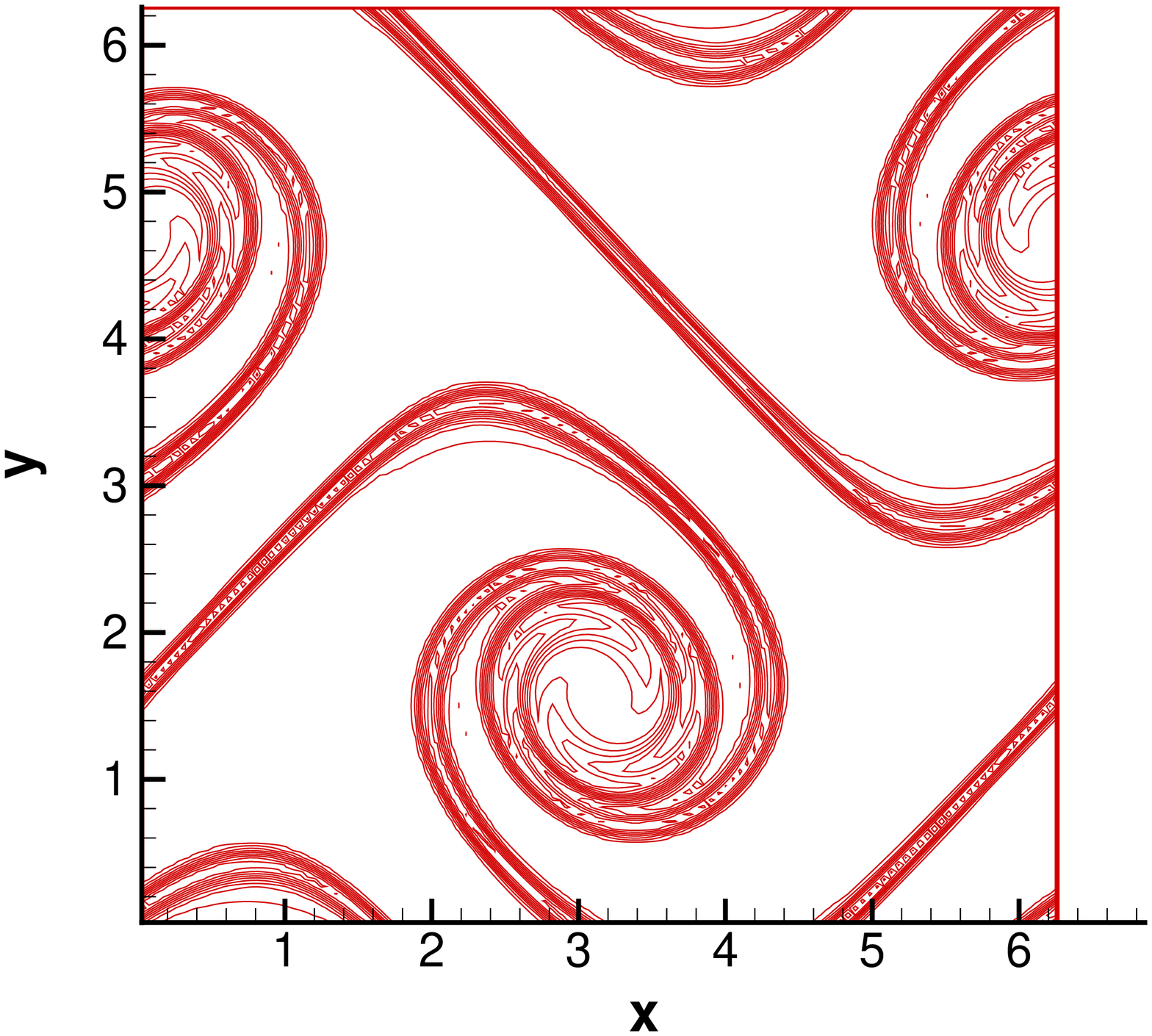}}
\caption{Numerical results for Example \ref{Eg:12} at $T=8$,
$64\times 64$(left),$128\times 128$(right).}\label{Fig:14}.
\end{figure}
\end{exa}

\begin{exa}\label{Eg:13}
(The vortex patch problem) Finally we solve the incompressible Euler equations (\ref{incom})
for the vortex patch problem on the domain $[0, 2\pi]\times [0, 2\pi]$ with the initial conditions
given by
\begin{equation}\label{ic:13}
\omega(x,y,0)=\left\{\begin{array}{ll}
-1, \quad & \textrm{$\frac{\pi}{2}\leq x \leq \frac{3\pi}{2}$,
              $\frac{\pi}{4}\leq x \leq \frac{3\pi}{4}$},\\
1, \quad & \textrm{$\frac{\pi}{2}\leq x \leq \frac{3\pi}{2}$,
              $\frac{5\pi}{4}\leq x \leq \frac{7\pi}{4}$},\\
0, \quad & \textrm{otherwise.}
\end{array}\right.
\end{equation}

In Fig. \ref{Fig:15}, we show the contour plots of vorticity and the
cuts along the diagonal at $T=5$ and $T=10$ for the FVCW scheme with
limiters. The results are also similar to those in
\cite{qiu2011conservative,zhang2010maximum,liu2000high}. The minimum
and maximum numerical values at $T=5$ with limiters are $-0.999995$
and $ 0.999995$, without limiters they are $ -1.000012$ and
$1.000012$. We omit the contour plots for the case without limiters
due to similarity.
\begin{figure}[htp]
\subfigure[T=5]{\includegraphics[width=0.5\textwidth]{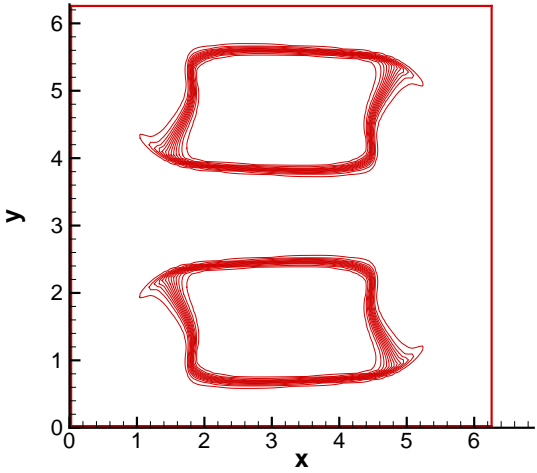}
\includegraphics[width=0.5\textwidth]{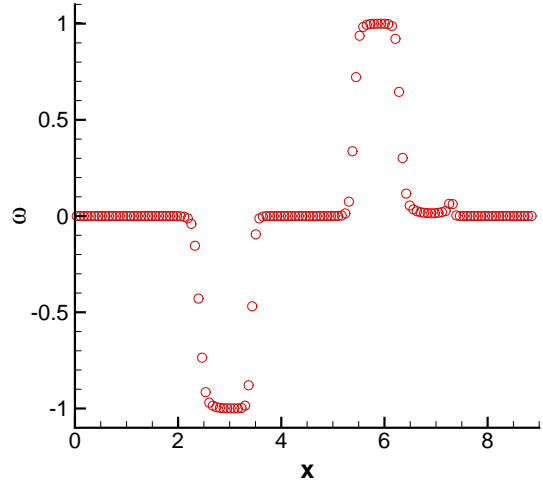}}
\subfigure[T=10]{\includegraphics[width=0.5\textwidth]{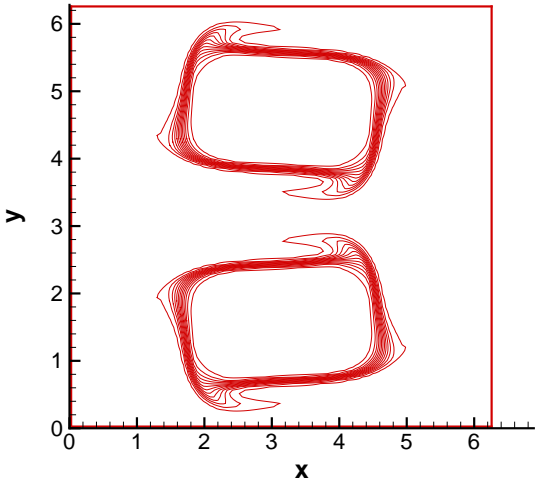}
\includegraphics[width=0.5\textwidth]{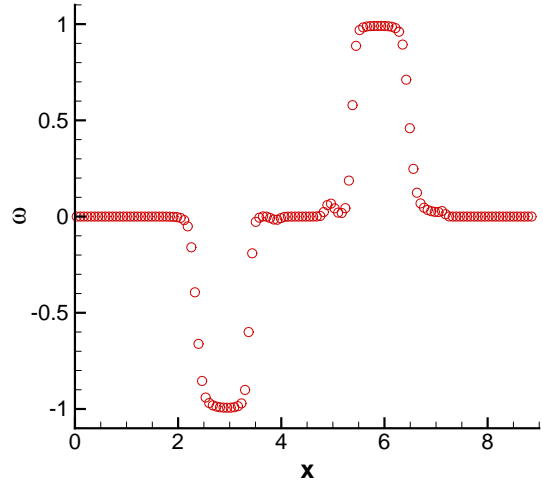}}
\caption{Numerical results of the FVCW scheme for Example
\ref{Eg:13}. Mesh $128\times 128$.}\label{Fig:15}.
\end{figure}
\end{exa}

%%%%%%%%%%%%%%%%%%%%%%%%%%%%%%%%%%%%%%%%%%%%%%ps

%%%%%%%%%%%%%%%%%%%%%%%%%%%%%%%%%%%%%%%%%%%%%%%
\section{Conclusions}
In this paper, we developed a maximum-principle-satisfying high order
finite volume compact WENO scheme. By applying a polynomial scaling limiter to the finite volume
compact WENO scheme at each stage of an explicit Runge-Kutta time method, the scheme satisfies
the strict
maximum principle under suitable CFL numbers without destroying high order of accuracy. Both one-dimensional and two-dimensional examples including incompressible flow problems are tested, the results showed that
the compact scheme has better resolution compared to the classical
finite volume non-compact WENO scheme. The application of our proposed method
on unstructured meshes is subject to our ongoing work.

\bigskip
\noindent
{\bf Acknowledgement.}
The work was partly supported by the Fundamental Research Funds for
the Central Universities (2010QNA39, 2010LKSX02). The third author
acknowledges the funding support of this research by the Fundamental
Research Funds for the Central Universities (2012QNB07).

\bibliographystyle{siam}
\bibliography{refer}

\end{document}